\documentclass[11pt]{amsart}
\usepackage{amssymb,amsmath,amsthm,amsfonts,mathrsfs}
\usepackage[hmargin=3cm,vmargin=3.5cm]{geometry}
\usepackage{hyperref}
\usepackage{color}
\usepackage{graphicx,psfrag}
\usepackage{amscd}  
\usepackage[shellescape]{gmp}
\usepackage{stmaryrd}  
\usepackage[all,2cell]{xy} \UseAllTwocells \SilentMatrices
\usepackage{tikz-cd}
\usepackage[dvips]{epsfig}
\usepackage{MnSymbol}

\newtheorem{theorem}{Theorem}[section]

\newtheorem{prop}[theorem]{Proposition}

\newtheorem{corollary}{Corollary}

\let\oldtocsection=\tocsection
\let\oldtocsubsection=\tocsubsection
\renewcommand{\tocsection}[2]{\hspace{0em}\oldtocsection{#1}{#2}}
\renewcommand{\tocsubsection}[2]{\hspace{1em}\oldtocsubsection{#1}{#2}}

\usepackage{chngcntr}
\counterwithin{figure}{subsection}

\title[Decorated 1D cobordisms]{Decorated one-dimensional cobordisms and tensor envelopes  of noncommutative recognizable  power series}

\author{Mikhail Khovanov}
 \address{Department of Mathematics, Columbia University, New York, NY 10027, USA}
 \email{\href{mailto:khovanov@math.columbia.edu}{khovanov@math.columbia.edu}}
 
\date{October 12, 2020}

\begin{document}

\def\R{\mathbb R}
\def\Q{\mathbb Q}
\def\Z{\mathbb Z}
\def\N{\mathbb N} 
\def\C{\mathbb C}
\def\S{\mathbb S}
\def\CP{\mathbb P}
\renewcommand\SS{\ensuremath{\mathbb{S}}}
\newcommand{\kllS}{\kk\llangle  S \rrangle} 
\newcommand{\kllSS}[1]{\kk\llangle  #1 \rrangle}
\newcommand{\klS}{\kk\langle S\rangle}  
\newcommand{\aver}{\mathrm{av}}  
\newcommand{\ophana}{\overline{\phantom{a}}}

\def\l{\lbrace}
\def\r{\rbrace}
\def\o{\otimes}
\def\lra{\longrightarrow}
\def\Ext{\mathrm{Ext}}
\def\mc{\mathcal}
\def\mf{\mathfrak} 
\def\mcC{\mathcal{C}}
\def\mcI{\mathcal{I}}
\def\uFr{\underline{\mathrm{Fr}}}

\def\bbb{\mathbb{B}}
\def\ovb{\overline{b}}
\def\tr{{\sf tr}} 
\def\det{{\sf det }} 
\def\tral{\tr_{\alpha}}
\def\one{\mathbf{1}}   

\def\lra{\longrightarrow}
\def\kk{\mathbf{k}}  
\def\gdim{\mathrm{gdim}}  
\def\rk{\mathrm{rk}}
\def\undep{\underline{\epsilon}}

\def\CCC{\mathcal{C}} 
\def\kCCC{\kk\CCC}  
\def\vCCC{\mathcal{V}\CCC}  
\def\sCCC{\mathcal{S}\CCC} 
\def\dCCC{\mathcal{D}\CCC} 
\def\udCCC{\underline{\dCCC}}  

\def\CCI{\widetilde{\mathcal{C}}} 
\def\kCCI{\kk\CCI}  
\def\vCCI{\mathcal{V}\CCI}  
\def\sCCI{\mathcal{S}\CCI} 
\def\dCCI{\mathcal{D}\CCI} 
\def\udCCI{\underline{\dCCI}}  

\def\Rec{\mathrm{Rec}} 

\def\Cob{\mathrm{Cob}} 
\def\Cobtwo{\Cob_2}   
\def\Kob{\mathrm{Kob}}
\def\Kar{\mathrm{Kar}}   

\newcommand{\alphai}{\alpha^{\bullet}}

\newcommand{\alphac}{\alpha^{\circ}}  
\newcommand{\alphap}{(\alphai,\alphac)} 

\newcommand{\oplusop}[1]{{\mathop{\oplus}\limits_{#1}}}
\newcommand{\ang}[1]{\langle #1 \rangle } 

\newcommand{\mcA}{{\mathcal A}}

\newcommand{\undn}{\mathbf{n}}
\newcommand{\cob}{\mathrm{cob}} 

\newcommand{\Hom}{\mbox{Hom}}
\newcommand{\id}{\mbox{id}}
\newcommand{\Id}{\mbox{Id}}
\newcommand{\End}{\mathrm{End}}

\begin{abstract} The paper explores the relation  between noncommutative power series and topological theories of  one-dimensional cobordisms decorated by labelled zero-dimensional submanifolds. These topological theories  give  rise to several types of tensor envelopes of noncommutative recognizable power series, including the categories built from the syntactic algebra and syntactic ideals of the  series and the analogue of the 
Deligne category. 
\end{abstract}

\maketitle
\tableofcontents

%
%

\section{Introduction}

In the \emph{universal construction} approach to low-dimensional topological theories~\cite{BHMV,Kh1,RW1} one starts with an evaluation of closed $n$-dimensional objects $M$ taking values in a  ground commutative ring or a field and then defines  state spaces $A(N)$ for $(n-1)$-dimensional objects $N$ via the bilinear pairing on $n$-dimensional objects $M$ with a given boundary, $\partial M \cong N$, by coupling two such objects $M_1,M_2$ along the  boundary and evaluating  the resulting  closed object $M_1 \cup_N M_2$. The $n$-dimensional objects  may be manifolds, manifolds with  decorations, embedded  manifolds or foams, or one  of many other variations of these examples. The universal pairing theory of Freedman, Kitaev, Nayak, Slingerland, Walker and Wang~\cite{FKNSWW}, further developed by Calegari, Freedman, Walker and others~\cite{CFW,W}, is  closely related to  the universal  construction. Some other examples of the  universal construction  for $n=2$ were recently considered in~\cite{Kh2,KS3,KQR}. Vector spaces or modules $A(N)$  that one  assigns to $(n-1)$-dimensional  objects in universal constructions usually do not satisfy  the Atiyah tensor product axiom $A(N_1\sqcup N_2)\cong A(N_1)\otimes A(N_2)$, see~\cite{A}. Instead, there are maps 
\[
\begin{tikzcd}
A(N_1)\otimes A(N_2) \arrow[r,""] &  
A(N_1\sqcup N_2)
\arrow[shift right, swap]{l}{}  
\end{tikzcd}  , 
\]
which one can  think of as a sort of a lax  tensor structure. 
 
 \vspace{0.1in} 
 
In this note  we  explain that  the universal construction approach is interesting even  in dimension one. Studying  the  universal construction for one-manifolds decorated  by dots labelled by elements of a finite set $S$, we recover the notion of noncommutative recognizable (equivalently, rational) power series in the alphabet $S$ as developed by Sch\"utzenberger~\cite{Sch}, Fliess~\cite{F}, Eilenberg~\cite{E1,E2}, Conway~\cite{Co}, Reutenauer, Carlyle  and  Paz, and others. A full set of references and introductions to  this theory can be found in the textbooks by Berstel and Reutenauer~\cite{BR2}, Salomaa and Soittola~\cite{SS},  \'Esik and  Kuich~\cite{EK}, Kuich and Salomaa~\cite{KuS}, also see~\cite{Sa2,HMS}. For  short introductions to noncommutative rational power series we refer  to Reutenauer~\cite{Re3,Re4,Re5}.

Theory of noncommutative recognizable power series has its roots in the theory of rational languages and finite state automata~\cite{Co,E1,E2,EK}, and can be viewed as a linearization of the latter~\cite{BR2,HMS}.
We briefly review the basics of noncommutative rational (recognizable) power series in Section~\ref{subset_nc_ser} and Proposition~\ref{prop_many} stated there. 
Part of the motivation for this theory comes from an earlier theorem of Kleene that rational languages are precisely those recognizable by FSA (finite state automata). A language $L$ is a subset of $S^{\ast}$ (the set of  words  in the letters  of the alphabet $S$)  and gives  rise to series $\alpha(L)$ with the coefficient of $w$ one if  $w\in  L$ and  zero  otherwise. Coefficients of $\alpha(L)$ belong
to the Boolean semiring $\bbb=\{0,1\}$ with $1+1= 1$. Kleene's theorem and  the theory of  finite state automata  can be found in many textbooks on the  field, see for instance~\cite{Ch} and foundational work of Conway~\cite{Co} and Eilenberg~\cite{E1,E2}. Alternatively, we refer to Underwood~\cite[Chapter 2]{U} for a brief introduction to finite state automata, regular languages, and  their relation to bialgebras.

\vspace{0.1in}

We consider various flavours of the category of $S$-decorated one-dimensional cobordisms. 
$S$-labelled dots placed along a one-dimensional cobordism  can  also be thought of as codimension one \emph{defects} on it. 

\vspace{0.1in}

In the first example, the category $\CCC$ of oriented one-dimensional  cobordisms with labels from $S$ is considered  in Section~\ref{sec_decorated}. We work over a ground field $\kk$ for simplicity, but  the construction extends to an arbitrary commutative ring $R$ and at least parts of it  extend to commutative semirings, mirroring the theory of recognizable power series over semirings.

To build an evaluation one needs a number (an element of the ground field $\kk$) associated to each circle carrying a collection of $S$-labelled dots. This collection is determined by a  finite  sequence $w$ of  elements of $S$ up to a cyclic order. Consequently, to build  various evaluation categories, we need to assign a number $\alpha(w)\in\kk$ to each  such sequence or word $w\in S^{\ast}$, subject to the  condition $\alpha(uv)=\alpha(vu)$ for all words $u,v$.

An evaluation of this  type  is  encapsulated by a formal expression 
\begin{equation}\label{eq_Z_rat}
    Z_{\alpha} = \sum_{w\in S^{\ast}}\alpha(w) \, w, \ \ 
     \alpha=\{\alpha(w)\}_{w\in S^{\ast}}, \ \ 
\end{equation} 
known as a noncommutative power series $Z_{\alpha}$, an element of the vector space $\kllS$ dual to the  free associative algebra $\klS$ generated by  elements  of  $S$: 
\begin{equation*}
    \kllS \ :=  \ \klS^{\ast}= \Hom_{\kk}(\klS,\kk). 
\end{equation*}

Recognizable noncommutative  power  series are singled out by  the condition that their \emph{syntactic algebra} $A_{\alpha}$,  see Section~\ref{subset_nc_ser}, is finite-dimensional. The syntactic algebra~\cite{Re1} is  the quotient of  $\klS$ by the largest  two-sided ideal $I_{\alpha}$ of  $\klS$ that lies in the hyperplane $\ker(\alpha)$,  when  $\alpha$  is considered as a linear map  $\klS\lra \kk$.

Property $\alpha(uv)=\alpha(vu)$ for all words $u,v\in S^{\ast}$ describes a particular type of  series that we refer to as \emph{symmetric} series. Reutenauer~\cite{Re1} calls such series \emph{central}. 
\vspace{0.1in} 

In Section~\ref{sec_decorated} we show  that for  recognizable symmetric series  $\alpha$ there is a satisfactory theory of tensor envelopes~\cite{Kn}, that  is, tensor  categories associated to $\alpha$, that mirrors the theory  of  the  Deligne  categories associated to symmetric groups and of negligible quotients of these categories~\cite{D,CO,EGNO}. Similar theories have recently  been introduced for evaluations of two-dimensional cobordisms in~\cite{Kh2,KS3}, two-dimensional cobordisms with corners~\cite{KQR}, and two-dimensional cobordisms with dots (codimension two defects)~\cite{KKO}. One can also compare our construction with tensor envelopes of the "one-sided inverse" algebras and  Leavitt path  algebras considered in~\cite{KT} for categorifications of rings of fractions and with the diagrammatic categorification of the  polynomial  ring in~\cite{KS1}. 

There  are several categories and functors between them associated to rational symmetric  noncommutative series  $\alpha$, defined throughout Section~\ref{subsec_from} and summarized in Section~\ref{subsec_summ} and diagram (\ref{eq_seq_cd_1}) there. Various skein and quotient categories that one  obtains extend the notion of the  syntactic algebra $A_{\alpha}$ of $\alpha$ (the quotient of noncommutative polynomials $\klS$  by the largest two-sided ideal contained in  $\ker\alpha$,  see above) and can be thought of as forming various tensor and  Karoubi  closures of the  latter. The theory of syntactic algebras of  noncommutative recognizable (or rational) power series was introduced and  developed by Reutenauer~\cite{Re1}. Syntactic algebra $A_{\alpha}$ appears as the endomorphism algebra of the generating object $(+)$ in several categories associated to $\alpha$. 

\vspace{0.1in}

In Section~\ref{sec_floating_cobs} we go beyond the restriction that  noncommutative power series be symmetric by enlarging  our category of  cobordisms. We consider  category $\CCI$  of  $S$-decorated cobordisms $M$ that  may have  endpoints strictly  "inside" the cobordism, that it, not on the top or bottom boundary $\partial_1 M$ and $\partial_0 M$. We call these \emph{inner}  or  \emph{floating} endpoints. Such  floating endpoints appear in diagrammatical calculi in~\cite{KS1,KT}, for instance. Cobordisms of this type between empty 0-manifolds (closed or floating cobordisms) have  connected components that are  either $S$-decorated oriented intervals or circles. A multiplicative evaluation on  such  cobordisms assigns an  element $\alphai(w)\in \kk$ to an oriented  interval with word $w$ written along it via labelled dots, and an element $\alphac(v)\in \kk$ to an oriented circle, with word $v$, well-defined up to cyclic rotation, written along it. 

Consequently, the analogue of  noncommutative power series in this  case is a pair
\begin{equation}\label{eq_a_pair_2}
\alpha=\alphap
\end{equation}
where $\alphai$ is a noncommutative power series and $\alphac$ is  a symmetric noncommutative power series. There does not  have  to be any relation between $\alphai$ and $\alphac$. 

Pair $\alpha$ of series as above allows  to evaluate $S$-decorated floating intervals (via $\alphai$) and floating circles (via $\alphac$). In  Section~\ref{sec_floating_cobs} 
to the evaluation data  $\alpha=\alphap$ we assign several tensor categories similar to those for the symmetric series. The resulting categories have the best behaviour when both $\alphai$  and $\alphac$ are recognizable series, and we specialize to this  case early. We  say that $\alpha$  is  recognizable if both  $\alphai$ and  $\alphac$  are  recognizable. 

\vspace{0.1in}

We follow the path familiar from Section~\ref{sec_decorated} and papers~\cite{KS3,KQR} and  assign several categories and  functors between them to each  recognizable pair $\alpha$, including the following categories:
\begin{itemize}
    \item  The category $\vCCI_{\alpha}$ of viewable cobordisms, where any closed (floating) component is reduced  via evaluation $\alpha$. 
    \item The skein category $\sCCI_{\alpha}$, where, additionally, elements of two-sided and one-sided syntactic ideals $I_{\alpha}$ and $I_{\alphai}^{\ell},   I_{\alphai}^{r}$ evaluate to zero when placed in the middle of the strand or by its floating endpoint, respectively. 
    \item The category $\CCI_{\alpha}$, the quotient  of  either $\vCCI_{\alpha}$ or  $\sCCI_{\alpha}$  by the ideal of negligible morphisms. 
    \item Additive Karoubi closure $\dCCI_{\alpha}$ of $\sCCI_{\alpha}$, which is the analogue of the Deligne category. 
    \item Additive Karoubi closure  $\udCCI_{\,\alpha}$  of $\CCI_{\alpha}$, equivalent to the  quotient of $\dCCI_{\alpha}$ by the ideal of  negligible morphisms. 
\end{itemize}

These four  categories have finite-dimensional hom spaces (again, assuming  $\alpha$ is recognizable), see diagram (\ref{eq_seq_cd_3}) and Section~\ref{sec_floating_cobs}. 
They can be thought of as various tensor envelopes of $\alpha$ and the syntactic algebra $A_{\alpha}$. 

\vspace{0.1in}

Categories built out of a single symmetric recognizable series in Section~\ref{sec_decorated} can be considered a special case of this construction, given by setting the first series $\alphai$ to zero. Setting the second  series $\alphac$ to zero, instead, results in another specialization of the theory, with all decorated circles evaluating to zero, while decorated intervals evaluating to coefficients of $\alphai$, see the remark at the end of Section~\ref{sec_floating_cobs}.

\vspace{0.1in}

In  this paper we use \emph{rational} and \emph{recognizable} interchangeably to refer to  noncommutative power series over a field with the syntactic ideal of finite codimension. Coincidence of rational and recognizable power series with coefficients in an arbitrary semiring is  a result of Sch\"utzenberger~\cite{Sch}, see also~\cite{BR2,SS,EK,KuS} for more details and references. For more general monoids, beyond  the  free monoid  on a finite set $S$, the sets of recognizable and rational series may differ, see~\cite{DG,Sa1} and references therein. The difference between rational and recognizable series is also visible in examples in~\cite{KQR}, where a recognizable series in two or more commuting variables needs to be rational with  denominators restricted to polynomials in single generating variables. 

\vspace{0.1in} 

The theory of  recognizable noncommutative power series makes sense over non-commutative semirings~\cite{BR2,SS}. 
One can look to generalize the theory  of tensor envelopes  of such series from  series over a field or a commutative ring to series over a semiring. Note that closed  cobordisms would then evaluate to elements of the ground semiring. Components of a closed cobordism "commute", in the  sense of  sliding past  each other, as elements of the commutative monoid of endomorphisms of  the unit object of the tensor category of cobordisms,  the empty zero-manifold. For this reason, it is natural to restrict to commutative semirings in this fuller extension of the
theory of tensor envelopes of noncommutative power series. We do not consider the general case of a ground commutative semiring $K$ in this paper, though, limiting ourselves to a ground field, but it may be  interesting   to develop.
The case when $K=\bbb$ is the boolean semiring, gives, in particular, the notion of tensor envelopes of a rational language $L$ or, equivalently, tensor  envelopes of a finite state  automaton. To get the  definition, run the constructions of Section~\ref{sec_floating_cobs} with $\bbb$ in  place of field $\kk$ and the pair $\alpha=(\alphai,0)$ of series with the zero symmetric series $\alphac =0$ and $\alphai$ the series of a regular language $L$. 
To test whether this notion is useful, one may study examples of quotients $\CCI_{\alpha}$  of the skein category $\sCCI_{\alpha}$ for such $\alpha$. 

\vspace{0.1in}

In the follow-up paper, we will consider one-dimensional cobordisms with more general  decorations, by edges and vertices of  an oriented graph (or a quiver) $\Gamma$. The graph $\Gamma$ may be finite or infinite. 
Dots on  a cobordisms are labelled  by oriented edges of $\Gamma$. Intervals of the cobordisms separated by dots along a connected component are labelled by vertices of $\Gamma$. A dot labelled by an edge $s:a\to b$ is surrounded by intervals labelled by vertices $a$ and  $b$, respectively, in the order that matches the  orientation of the corresponding connected component. Such decorations  are possible for both interval and circle connected components of a cobordism. There are suitable monoidal  categories $\CCC(\Gamma)$ and $\CCI(\Gamma)$ of $\Gamma$-decorated cobordisms generalizing  categories $\CCC$ and  $\CCI$ in this paper. Cobordisms with floating endpoints are allowed in $\CCI(\Gamma)$ but   not in $\CCC(\Gamma)$. Objects of $\CCC(\Gamma)$ and $\CCI(\Gamma)$ are  finite  sequences of vertices of $\Gamma$ or, equivalently,  finite sequences  of objects of category  $S(\Gamma)$, see next. 

To $\Gamma$ one assigns the small category $S(\Gamma)$ of paths in $\Gamma$, with vertices of $\Gamma$ being the  objects of $S(\Gamma)$ and paths in $\Gamma$ -- morphisms, with concatenation of paths as the composition. Traveling  along a connected component of a $\Gamma$-decorated cobordism one encounters a path in  $\Gamma$, that is,  a morphism  in $S(\Gamma)$. If the  component is a circle, the path, in addition,  must be  closed, that is, start  and end at the same vertex of $\Gamma$.

An  evaluation $\alpha$, in the case of $\CCI(\Gamma)$, where floating endpoints are allowed,  consists of two maps: 
\begin{itemize}
    \item Map $\alphai$ from the set of  morphisms  in $S(\Gamma)$ (paths in $\Gamma$)  to the ground  field  $\kk$ or, more generally, a commutative ring or a semiring,  
    \item Map $\alphac$ from the set  of \emph{circular  morphisms}, that  is, closed paths in $\Gamma$ without a choice of  the basepoint to $\kk$. 
\end{itemize}
The pair $\alpha=(\alphai,\alphac)$ is  the  analogue of the pair  in (\ref{eq_a_pair_2}), generalizing the special case considered  in this paper where $\Gamma$ has a single vertex and oriented loops from the  vertex to itself are enumerated by elements of $S$.

One can  then define the analogues of all the categories in the diagram (\ref{eq_seq_cd_3}), including $\vCCI_{\alpha}$, $\sCCI_{\alpha}$,  $\CCI_{\alpha}$, in this case. In particular, the  category $\CCI(\Gamma)_{\alpha}$ is the quotient of
the $\kk$-linearization $\kk\CCI(\Gamma)$ by the two-sided ideal of negligible morphisms, defined via the trace given by evaluation $\alpha$. 

The pair $\alpha$ is called \emph{recognizable} or  \emph{locally-recognizable} (when $\Gamma$ is infinite) if the "gligible quotient" category $\CCI(\Gamma)_{\alpha}$  has finite-dimensional hom spaces.  Note that the boundary points and floating  endpoints of $\Gamma$-decorated cobordisms are labelled by  objects of $S(\Gamma)$, that is,  by vertices of $\Gamma$, the label inherited  from the label of  the adjacent edge of the cobordism. 

These constructions can be further generalized to cobordisms between finite sets of boundary points given  by graphs  and $\Gamma$-decorated  graphs rather then by $\Gamma$-decorated one-manifolds. Cobordisms given by graphs  can still be viewed as one-dimensional cobordisms between zero-dimensional  objects (finite sets of points, possibly decorated by vertices of  $\Gamma$ and orientations, as necessary). 

\vspace{0.1in}

 A natural open problem is to extend the universal construction, for decorated cobordisms (or cobordisms with defects), beyond dimension one. Parts of this extension are visible in
 \begin{itemize} 
 \item \cite{KKO}, where two-dimensional cobordisms  are decorated by dots labelled by elements of a commutative monoid  or a commutative algebra, with non-trivial interactions between these dots and topology of cobordisms coming from the handle cobordism equated to a nontrivial element  of the monoid or algebra. 
 \item \cite{KQR}, where side boundaries of two-dimensional cobordisms  with corners may be colored by elements of a finite set. 
 \item Foam theory~\cite{Kh1,B,EST,RWd,RW1}, see more references in~\cite{KK}, where rather particular evaluations of two-dimensional decorated CW-complexes with  generic singularities embedded  in $\R^3$ (foams) are used as an intermediate step to  build homology theories of links that categorify various one-variable specializations of the HOMFLYPT polynomial. Soergel bimodules, singular Soergel  bimodules, and some other structures in representation theory admit a foam description as well~\cite{Wd,RW2,KRW}.  
 \item Evaluation theory for two-dimensional cobordisms and evaluations of overlapping foams~\cite{Kh2}, pointing towards further connections to arithmetic topology, representation theory, and the Heegaard-Floer theory. 
 \item \cite{KL}, which considers  evaluations in the   two-dimensional planar case with one-dimensional defects. 
 \end{itemize} 
 Freedman et. al~\cite{FKNSWW} mention possible decorations on low-dimensional cobordisms for their universal pairings.
 
 \vspace{0.1in} 
 
 Studying evaluations not just for $n$-manifolds but for decorated $n$-manifolds, $n$-manifolds and their foam analogues embedded  in $\R^{n+1}$, and other such refinements should ease one's way  into understanding recognizable evaluations in dimension $n+1$. This program makes sense at least in dimensions $n=1,2,3$. 
 
\vspace{0.1in} 

{\bf Acknowledgments.} The author is grateful  to  Kirill  Bogdanov, Mee Seong Im, and Vladimir Retakh for interesting discussions and to Victor Shuvalov for help with creating the figures for the paper. The author was partially supported by the NSF grant  DMS-1807425 while working on this paper. 

%
%

\section{Decorated one-dimensional cobordisms and their  evaluations}\label{sec_decorated}


\subsection{Categories \texorpdfstring{$\CCC$}{C} and \texorpdfstring{$\CCC'$}{C'}}

Fix a finite set $S$ of cardinality $r\ge 0$, which  we  often write  as $S=\{s_1,s_2,\dots, s_r\}$. 

\begin{figure} 
\begin{center}
\includegraphics[scale=1.0]{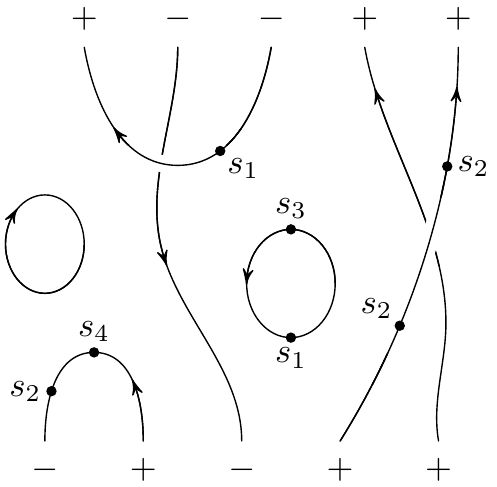}
\caption{A morphism from $(-+-++)$ to $(+--++)$ in $\CCC$.  It  has two closed (or floating) components, one  undecorated, the  other decorated by  $s_1s_3$. The two U-turns are decorated by $s_2s_4$ and $s_1$,  respectively. There  are three through arcs: two undecorated, one decorated  by $s_2s_2$. Whether  a crossing is over- or under-crossing is irrelevant.  }
\label{fig1_1}
\end{center}
\end{figure}

Consider  the category $\CCC=\CCC_S$ of $S$-decorated  compact oriented one-dimensional cobordisms. Its objects are oriented one-dimensional manifolds $N$, that is, finite  sets with a sign assignment $+$  or $-$  to each element (signed  finite sets). A morphism from $N_0$ to $N_1$  is  an oriented  one-dimensional manifold  $M$ decorated by finitely many dots  labelled by  elements of $S$,  with $\partial M = N_1 \sqcup (- N_0)$, see Figure~\ref{fig1_1} for  an example, 
which also sets the orientation convention for the boundary. 

Dots can move  along a connected component of $M$ where  they are  placed but without crossing through other dots or moving  to  a boundary point. Two morphisms are equal if they are diffeomorphic  rel boundary and  keeping track of dots  and  their labels. 

Each component $c$ of $M$  is  either an  oriented  circle or an oriented interval. Going along  $c$ in the  direction  of its orientation, one can read off  the labels of marked points. When $c$ is an interval, the sequence of labels is  an  invariant  of $c$. When $c$ is a circle, the sequence of labels is  an invariant up to a cyclic rotation of the sequence. A component  may carry no dots; the corresponding sequence is empty then. 

Composition of morphisms is given by their  concatenation.  

\vspace{0.1in} 

To reduce to fewer  objects, we take the objects to be sequences of signs $\undep=(\epsilon_1,\dots, \epsilon_n)$, $\epsilon_i\in \{+,-\}$. To $\undep$ we associate an ordered signed zero-manifold with one point for each term in the sequence, with the signs given by  $\epsilon_1,\dots, \epsilon_n$. We may alternatively  write $\epsilon_i=1$ or $-1$ instead of $+$ or $-$. 

Permutation cobordisms show that permuting signs in a sequence  $\undep$ leads to an isomorphic object, and that isomorphism classes of objects are parametrized by pairs of non-negative  integers $\undn=(n_0,n_1)$, counting  the number  of plus and minus signs. 

When restricting to  a skeleton category (one object for each isomorphism class), we thus reduce objects to pairs  $\undn=(n_0,n_1)$,  where $n_0$ is  the number of plus points and  $n_1$ is the number  of minus points. In the sequence  of signs that $\undn$ represents, we put plus signs  first, and  can also  write  $\undn=(+^{n_0}-^{n_1})$. 

\vspace{0.1in} 

Denote  by $\CCC=\CCC_S$ the category of $S$-decorated cobordisms  with  objects---finite sign sequences $\undep$ as above. The skeleton category of $S$-decorated cobordisms, with objects $\undn$, is denoted  $\CCC'$.
Category $\CCC$ is slighter  larger than the equivalent   category $\CCC'$. 

These  categories  are rigid  symmetric tensor, with  the  tensor  product in $\CCC$ given on morphisms by  placing their diagrams next  to  each other. On objects,  the tensor  product  is the concatenation of sequences.
In $\CCC'$ when forming the tensor product of morphisms, we  group plus points together  and  minus points  together.  Tensor product on objects  is  given by $(n_0,n_1)\otimes (m_0,m_1)=(n_0+m_0,n_1+m_1)$.  

\begin{figure} 
\begin{center}
\includegraphics[scale=1.0]{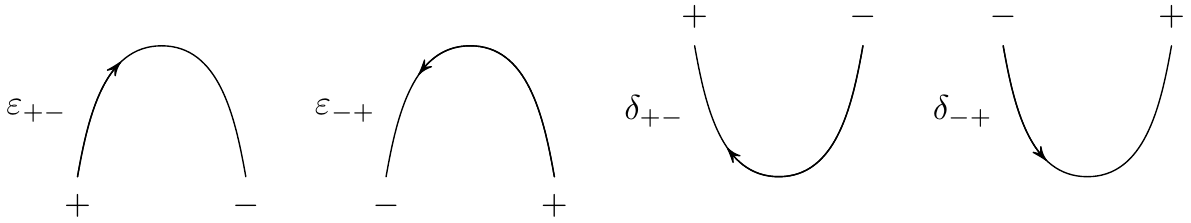}
\caption{Duality morphisms for $+$ and  $-$ in $\CCC$.  }
\label{fig1_2}
\end{center}
\end{figure}

In both categories $\CCC$ and  $\CCC'$ object $(+)$ has object $(-)$ as its dual, with the duality morphisms in  $\CCC$ shown in Figure~\ref{fig1_2}. 

\begin{figure} 
\begin{center}
\includegraphics[scale=1.0]{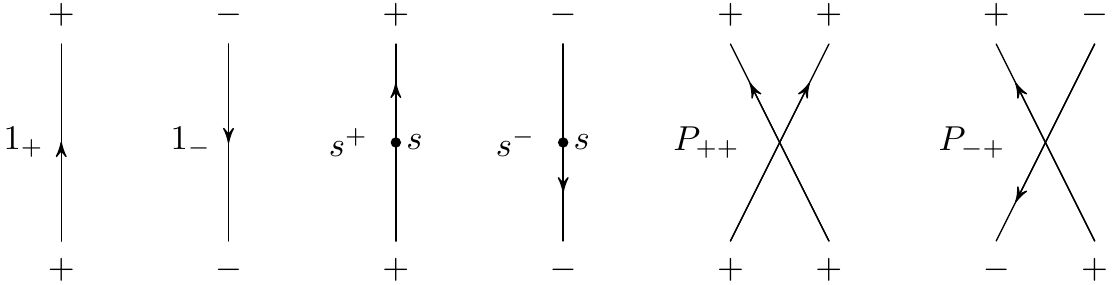}
\caption{Generating morphisms in $\CCC$: identity $1_+$ of $(+)$, identity $1_-$ of $(-)$, morphism $s^+:(+)\lra (+)$ for $s\in S$, morphism $s^-: (-)\lra (-),$  permutation  morphismm $P_{++}$. Other permutation morphisms,  such as $P_{-+}$, can be obtained as compositions of these morphisms and  those  in Figure~\ref{fig1_2}.}
\label{fig1_3}
\end{center}
\end{figure}

The  empty sequence  $\emptyset$ is the unit object $\one$ of the tensor  category  $\CCC$. The pair  $\one:=(0,0)$ is the unit object  of $\CCC'$.
Generating morphisms in $\CCC'$ are shown in Figures~\ref{fig1_2} and~\ref{fig1_3}.

\begin{figure} 
\begin{center}
\includegraphics[scale=1.0]{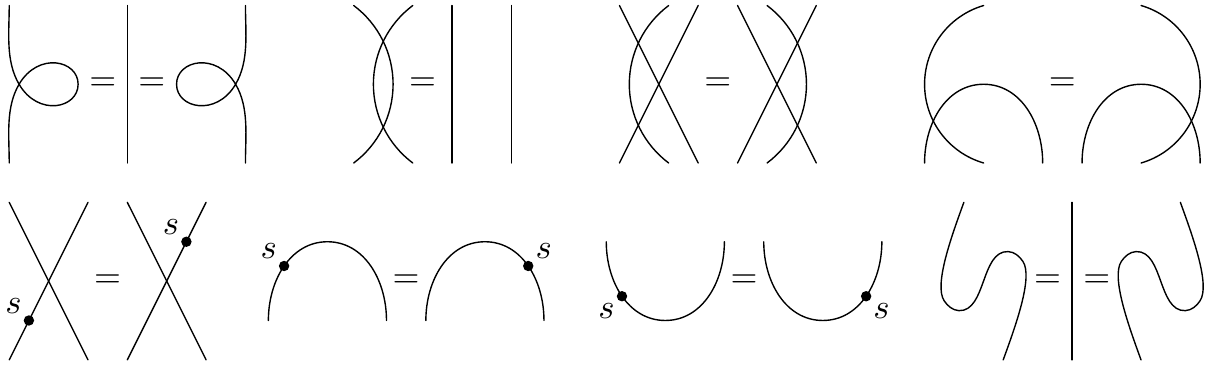}
\caption{Some relations in $\CCC$. They hold for any choice  of  orientations. Orientation of the  LHS of  each equation  determines the orientation of the RHS and vice versa. A set of these relations with some restrictions on orientations can be taken for a defining set of  relations. Commutativity relations on generators in  horizontally separated regions are not  shown, since they are built into the axioms of a tensor category.}
\label{fig1_4}
\end{center}
\end{figure}

Some defining relations in $\CCC$ are  shown  in Figure~\ref{fig1_4}. We do not list a full  set of defining relations and  will not need it. These relations can be hidden in the definition of $\CCC$, where morphisms are declared equal if the corresponding  decorated cobordisms are diffeomorphic  rel boundary. 

A morphism $M$ from $\epsilon$ to $\epsilon'$ in $\CCC$  consists of  some number of oriented  circles and oriented intervals. Boundaries of oriented  intervals and their orientations match entries  of $\epsilon$  and  $\epsilon'$ in pairs. 

Denote by  $|\epsilon|$ the difference of the  number of  plus and minus  signs  in  $\epsilon$ and call it the weight of the sequence. For instance, $|(+--++++-)|=5-3=2$. 
A morphism  from $\epsilon$  to $\epsilon'$ exists iff the two sequences  have the same weight, $|\epsilon|=|\epsilon'|$. The weight is  additive under  the tensor product  of objects (concatenation  of sequences). Denote  by  $||\epsilon ||$ the  length of the  sequence $\epsilon$. 

Connected components $c$  of a morphism $M$ are circles and intervals (arcs). Circles are closed components, also called \emph{floating} components. Arcs have boundary and, borrowing terminology from~\cite{KS2}, separate  into \emph{\it U-turns}  and  \emph{through} arcs. A U-turn has both endpoints on the same side of a morphism (either on the source one-manifold or on the target, while a \emph{through} arc has  one  endpoint  on  the source and  one on the target, also see types (1)-(3) of components in Figure~\ref{fig4_3}. The endpoints of a U-turn have opposite signs, while the  endpoints of a through  arc carry  the  same sign, see Figures~\ref{fig1_2},~\ref{fig1_3}. 

By analogy with~\cite{KS3,KQR}, we can also call arcs \emph{viewable} or \emph{visible} components, since they  have endpoints on the boundary of the cobordism (either top or bottom or both), and  call circles \emph{floating} components~\cite{KS1,KQR}, since they are disjoint from the boundary of the cobordism. 

\vspace{0.1in} 

Denote by $S^{\ast}$ the set of finite sequences of elements of $S$, including  the empty sequence $\emptyset$ (also see Section~\ref{subset_nc_ser}).  Going along an arc $c$ of a cobordism $M$ gives us  a word $w(c)\in S^{\ast}$. Going along a circle $c$ in $x$ gives a word  $w(c)$ well  defined up  to a  cyclic rotation or \emph{conjugation}
of words, $w_1w_2\sim w_2w_1$. 

In this  paper we encounter sequences $\undep$  of  signs, which are objects of $\CCC$, and  sequences $w\in  S^{\ast}$, which are sequences  of labels encountered along connected components  of a cobordism,  the latter a morphism in $\CCC$.


\subsection{Noncommutative power series}\label{subset_nc_ser}
$\quad$
\vspace{0.1in} 

For simplicity we work over a ground  field $\kk$, although the theory of noncommutative power series and rational and recognizable series makes sense over an arbitrary semiring $R$, not necessarily  commutative~\cite{BR2,SS}. For definitive treatments we refer  the reader to books~\cite{BR2,SS} and to~\cite{Re3,Re4,Re5} for quick introductions  and reviews.  

Let $S^{\ast}=\emptyset  \sqcup S\sqcup S^2\sqcup \dots $ be the set of sequences $w=t_1\dots t_n$  of  elements  of a  finite set $S=\{s_1,\dots,  s_r\}$. We call elements  of $S$ \emph{letters}  and elements of  $S^{\ast}$ \emph{words} or \emph{sequences} in $S$.  The empty word $\emptyset$ is allowed. A noncommutative power series $\alpha$ over $\kk$ is any function 
\begin{equation}
    \alpha \ :  \ S^{\ast} \lra \kk, \ \ \alpha(w)\in \kk,  \ w\in  S^{\ast}.
\end{equation}
We formally write  this  series as 
\begin{equation}\label{eq_Z_rat_2}
    Z_{\alpha} = \sum_{w\in S^{\ast}}\alpha_w \, w, \ \ 
     \alpha=(\alpha_w)_{w\in S^{\ast}}, \ \ \alpha_{w}=\alpha(w) \in \kk , 
\end{equation} 
using either $\alpha_w$ or  $\alpha(w)$ to denote the value of  $\alpha$ on a noncommutative monomial or  word $w$. 
Denote by $\kllS$  the $\kk$-vector space of noncommutative power  series  and by $\klS$ the free noncommutative $\kk$-algebra on generators in $S$ (the  algebra of noncommutative polynomials). 

Given two series $\alpha,\beta$, their product is the series $\alpha\beta$ that on $w$ evaluates to
\begin{equation}
    \alpha\beta(w)=\sum_{w=w_1w_2}\alpha(w_1)\beta(w_2), 
\end{equation}
the  sum over all decompositions of $w$. There are $\ell(w)+1$ terms in the sum, where $\ell(w)$ is the length of $w$. This product turns $\kllS$ into a $\kk$-algebra, noncommutative if $S$ has  more than  one element. The  inclusion $\klS\subset \kllS$ is a ring homomorphism. 

\vspace{0.1in} 

We say that  series $\alpha\in \kllS$  is
\emph{recognizable} iff there is a homomorphism $\psi:\klS\lra \End(\kk^n)$ of the free algebra into the  algebra of $n\times n$ matrices, a  vector and a dual vector $\lambda,\mu^T\in \kk^n$ such that 
\begin{equation}\label{eq_recog}
    \alpha(w) = \mu \, \psi(w)\, \lambda 
\end{equation}
for all words $w$. That  is, the  number $\alpha(w)$ is the product of the $1\times n$ matrix $\mu$, $n\times n$ matrix $\psi(w)$ and $n\times 1$ matrix $\lambda$.
Denote by $\kllS^{rec}$ the set  of all recognizable series. 

\vspace{0.1in} 

Vector space $\kllS$ is a $\klS$-bimodule with  $f\otimes g\in \klS\otimes  \klS^{op}$ acting on $\alpha \in  \kllS$ by
\begin{equation*}
    (f\otimes g)(\alpha)(w) = \alpha(gwf) , \ \  w\in S^{\ast}.
\end{equation*}
We write $f\alpha g := (f\otimes g)(\alpha)$. 
This action gives rise to the left, right, and two-sided ideals $I^{\ell}_{\alpha}, I^{r}_{\alpha}$ and $I_{\alpha}$ in $\klS$: 
\begin{itemize}
    \item 
Left ideal $I^{\ell}_{\alpha}$  consists of  all $f\in   \klS$ such  that  $f\alpha=0,$ that  is, all $f$ such that $\alpha(wf)=0$ for any word $w\in  S^{\ast}$. It  is the largest left ideal contained in the hyperplane $\ker(\alpha)\subset \klS$.  
\item 
Right ideal $I^r_{\alpha}$   consists of  all $g\in   \klS$ such  that  $\alpha  g=0$,  that  is, all $g$ such  that $\alpha(gw)=0$ for all $w\in \S^{\ast}$. It  is the largest right ideal contained in $\ker(\alpha)$.  
\item 
Ideal $I_{\alpha}$ consists of all $f\in \klS$ such that $\alpha(wfv)=0$ for all $w,v\in S^{\ast}$. It is the largest two-sided ideal  contained  in the hyperplane $\ker(\alpha) $. 
\end{itemize}

Ideal $I^{\ell}_{\alpha}$ has finite codimension in $\klS$ iff the series $\alpha$  is recognizable. 
Given triple as in (\ref{eq_recog}), ideal $I^{\ell}_{\alpha}$ contains  the  finite codimension subspace $\{ x\in \klS| \psi(x)\lambda=0\}$. 
Vice versa, if $I^{\ell}_{\alpha}$ has finite codimension, it is straightforward to produce the data in (\ref{eq_recog}) by taking
$\kk^n\cong \klS/I^{\ell}_{\alpha}$, $\lambda=1$ and $\mu=\alpha$. 
Given a triple as in (5), $\ker\alpha$  contains two-sided ideal of finite  codimension $\{f\in \klS |  \psi(f) =0\in \End(\kk^n)\}.$ 
Vice versa,  if $I_{\alpha}$ has finite codimension, ideals $I^{\ell}_{\alpha}\supset  I_{\alpha}$  and  $I^r_{\alpha}\supset I_{\alpha}$ have finite codimension too. 

Consequently, if one of $I^{\ell}_{\alpha},I^r_{\alpha},I_{\alpha}$ have finite codimension in $\klS$, the other two have finite codimension as well.

Two-sided  ideal $I_{\alpha}$ of $\klS$ is  called the \emph{syntactic ideal} of $\alpha$. Denote by
\begin{equation}\label{eq_synt_alg}
    A_{\alpha}  \ := \ \klS/I_{\alpha}
\end{equation}
the quotient algebra, the \emph{syntactic algebra} of $\alpha$, see~\cite{Re1}. It is defined for any $\alpha,$ but we mostly restrict to considering it for recognizable  $\alpha$, when  $A_{\alpha}$  is finite-dimensional. We also call $I^{\ell}_{\alpha}$  and  $I^{r}_{\alpha}$ the \emph{left} and \emph{right syntactic ideals} of  $\alpha$. 

Algebra $\klS$ acts on $\kllS$ on the  left  and on the  right, and $\klS$-bimodule generated by $\alpha$ (a subbimodule  of $\kllS$) is naturally isomorphic  to the syntactic algebra $A_{\alpha}$, the latter equipped with $\klS$-bimodule structure via  left and  right  multiplications: 
\begin{equation}
    A_{\alpha} \cong \klS\otimes  \klS^{op}(\alpha).
\end{equation}

The  quotient   $\klS/I^{\ell}_{\alpha}$ is naturally a faithful
left $A_{\alpha}$-module via the left multiplication action. 
Likewise, $\klS/I^{r}_{\alpha}$ is a faithful right $A_{\alpha}$-module via the right multiplication action. 

We see that  $\alpha\in \kllS$  is recognizable iff the cyclic $\klS\otimes  \klS^{op}$-module generated  by  $\alpha$ in $\kllS$ is finite-dimensional, or, equivalently, 
\begin{equation*}
    \dim_{\kk} A_{\alpha} < \infty. 
\end{equation*}

The \emph{Hankel matrix} $M_{\alpha}$  of $\alpha$  is the infinite square matrix with rows and  columns enumerated by elements  of $S^{\ast}$ with the $(w_1,w_2)$-entry $\alpha(w_1w_2)$. 

\vspace{0.1in} 

Given any series $\alpha$ with $\alpha(\emptyset)=0$ (called \emph{proper series}), we can form the \emph{Kleene plus} series $\alpha^+$ as the formal sum
\begin{equation}
    \alpha^+ = \alpha  + \alpha^2+\dots,
\end{equation}
where $\alpha^n=\alpha\alpha\dots \alpha$ is the product of $n$ copies of $\alpha$. The term $\alpha^n$ evaluates to $0$ on any word of length less  than $n$. Consequently, a given word evaluates  nontrivially only on finitely many terms in the  sum, and $\alpha^+$ makes sense as an element of $\kllS$. The series $1+\alpha^+$ is the inverse of the series $1-\alpha$ in the ring $\kllS$. 

A series is \emph{finite} if it contains finitely many terms. Finite series are those in the ring of noncommutative polynomials $\klS\subset \kllS$. 

Denote by $\kllS^{rat}$ the smallest subset of series that 
\begin{itemize}
    \item Contains all finite series. 
    \item Closed under the product and finite $\kk$-linear combinations of series. 
    \item Contains $\alpha^+$ for any proper series $\alpha$ in the subset. 
\end{itemize}
Series in $\kllS^{rat}$ are called \emph{rational} series. 

\begin{prop} \label{prop_many} 
The  following properties of series $\alpha$ are equivalent. 
\begin{enumerate}
    \item $\alpha$ is rational.
    \item $\alpha$  is recognizable. 
    \item The  Hankel matrix $M_{\alpha}$ of $\alpha$ has finite rank.  
    \item The  syntactic ideal $I_{\alpha}$  has finite  codimension  in $\klS$. 
    \item The left ideal $I^{\ell}_{\alpha}$ has finite codimension in $\klS$. 
    \item The right ideal $I^{r}_{\alpha}$ has finite codimension in $\klS$. 
    \item $\alpha$  can be computed by a weighted finite automaton. 
\end{enumerate}
\end{prop}

Equivalence of (2), (4), (5), (6) is explained above. 

For a proof of all equivalences see Sections 1 and 2 of~\cite{BR2}, Salomaa-Soitola~\cite{SS}, or references there to  the original work of Sch\"utzenberger~\cite{Sch}, Fliess~\cite{F}, Eilenberg~\cite{E2} and others. Most of these equivalences hold in much greater generality than over a field, in many cases over an arbitrary semiring. The Hankel matrix of noncommutative series was introduced  by Fliess~\cite{F}.  

The notion of \emph{weighted finite automaton} linearizes the concept of finite state automaton and, over a field $\kk$,  is equivalent to the triple $(\lambda,\psi,\mu)$ as in (\ref{eq_recog}), see~\cite[Section 1.6]{BR2}, for  instance. 
$\square$

We have $\kllS^{rec}=\kllS^{rat}$, since rational and recognizable series coincide. 

Assume  that $\alpha$  is recognizable. 
The trace  form $\alpha$  on the  finite-dimensional algebra $A_{\alpha}$ has the  following nondegeneracy property:
\begin{equation}
\label{eq_trace1}
\mbox{for  any $a\in A_{\alpha}, a\not=0$ there are $b,c\in A_{\alpha}$ such that  $\alpha(bac)\not= 0$.} 
\end{equation}

This is a much weaker condition than the usual Frobenius  condition on a linear form $\beta$ on a finite-dimensional  algebra $B$:
\begin{equation}
\label{eq_trace2}
\mbox{for any $a\in B,$ $a\not=0$ there exist $b$ such  that $\beta(ab)\not=0$.}
\end{equation} 
In the latter case $\beta$ equips $B$ with the  structure of a Frobenius  algebra.

Given any finite-dimensional algebra $B$ with a linear form $\alpha: B\lra \kk$,  the condition  that  
\begin{equation}\label{eq_trace3}
\mbox{
for  any $a\in B, a\not=0$ there are $b,c\in B$ such that  $\alpha(bac)\not= 0$}
\end{equation}
is equivalent to the zero ideal $(0)$ being the only two-sided ideal in $\ker(\alpha)$.  Let  us  call a pair $(B,\alpha)$  with this property a \emph{syntactic pair}. A finite set of  generators $b_1, \dots, b_m$ of $B$ gives  rise to a surjective homomorphism  
\begin{equation}
    \rho: \klS\lra B,\ \ \rho(s_i)=b_i, \ \ S=\{s_1,\dots, s_m\}
\end{equation}
from the free algebra $\klS$ to $B$ and induced noncommutative power series in the set of variables  $S$, also denoted $\alpha$. This gives a bijection between recognizable power series in $S$ and isomorphism  classes of  syntactic pairs $(B,\alpha)$ as  above with a  choice of generators $(b_1,\dots, b_m)$ of $B$. 

\vspace{0.1in} 

An algebra is  called \emph{syntactic} if it admits a presentation (\ref{eq_synt_alg}) for some $S$  and $\alpha$. 

{\emph Examples:}
\begin{enumerate}
    \item Any Frobenius algebra $B$ with a non-degenerate  form $\beta$ gives a syntactic pair $(B,\beta)$. 
    \item Take the matrix algebra $B=\mbox{M}_n(\kk)$ and define $\alpha(x)=x_{1,1}$ to pick the first diagonal coefficient of the matrix $x$. 
    In this example the form $\alpha$ satisfies the  weaker property (\ref{eq_trace1}), so that $(B,\alpha)$ is a syntactic pair, but $\alpha$ is not a  Frobenius trace. Algebra $B$ is Frobenius for a different linear form on it (for example, for the usual trace on matrices).
    \item Take the path algebra $B$ of  the quiver  with two  vertices $0,1$  and the edge $(01)$ connecting  them, 
    with the  multiplication given by concatenation of paths:  $(0)(01)=(01),(01)(1)=(01)$, etc. Algebra $B$  has  a basis $\{(0),(1),(01)\}$. Take any linear  form  $\alpha$ with $\alpha((01))\not= 0$. Then  $(B,\alpha)$ is a syntactic algebra with  this linear form. It can be generated by two elements. $B$ is neither Frobenius nor quasi-Frobenius. 
    \item A finite-dimensional commutative algebra  is Frobenius iff it is syntactic.  
\end{enumerate}

See Reutenauer~\cite{Re1} and Perrin~\cite{Pe} for more results on syntactic algebras and the latter also for another brief introduction to the subject.

\vspace{0.1in}


\subsection{Evaluations and symmetric series} 

We say that series $\alpha \in  \kllS$ is \emph{symmetric} if $\alpha(w_1w_2)=\alpha(w_2w_1)$ for any $w_1,w_2\in S^{\ast}$. Transformation $w_1w_2\mapsto w_2w_1$ is also called \emph{conjugation}, so one  can say  that   $\alpha$ is conjugation  invariant. An evaluation $\alpha$ is symmetric iff it only depends on a sequence up to cyclic order. 

We use the word "symmetric"  to define such series, since the word "cyclic" is already taken, see~\cite{BR1,KaR,Re2} and~\cite[Section 12.2]{BR2}. A series $\alpha$ is called \emph{cyclic} if,  in addition to the conjugation invariance condition, it  satisfies $\alpha(w^n)=\alpha(w)$ for any non-empty $w$. Thus, a cyclic series is symmetric but most symmetric series are not cyclic. Reutenauer~\cite{Re1} uses \emph{central} instead of our \emph{symmetric}.  

Denote  the set of symmetric series by $\kllS^s$ and by $\kllS^{s,rec}$ the set of recognizable symmetric series.

A series $\alpha\in \kllS$ can be averaged out  to a series $\aver(\alpha)$ given by
\begin{equation}
    \aver(\alpha)(w) \ =  \ \sum_{uv=w,v\not=\emptyset} \alpha(vu), \ \ \mathrm{if} \ w\not= \emptyset,\ \ \aver(\alpha)(\emptyset)=\alpha(\emptyset).
\end{equation}
That  is,  take the sum over all possible ways to split   $w$  into the product  $uv$ and evaluate  $\alpha$ on $vu$. Series $\aver(\alpha)$  is symmetric. Only one of the two  degenerate splittings $\emptyset w$ and $w\emptyset$  is used to avoid having $\alpha(w)$ twice  in the sum.  

\begin{prop}   $\aver(\alpha)\in\kllS^{s,rec}$ if $\alpha\in \kllS^{rec}.$
\end{prop} 
In other words, averaging out a recognizable series  produces  a symmetric recognizable series. This result is proved in Rota~\cite{Ro}, see also~\cite{Re2}. It gives a large supply  of symmetric recognizable series. 
$\square$ 

\vspace{0.1in} 

Symmetric series with semisimple syntactic algebra $A_{\alpha}$  are studied in~\cite{Re1,Pe}.

\vspace{0.1in}


\subsection{Tensor envelopes  of  series \texorpdfstring{$\alpha$}{a}} \label{subsec_from}
\quad
\vspace{0.1in} 

{\it (1) Category  $\kk\CCC$.}
We fix a base field $\kk$ and form  the $\kk$-linearization $\kk\CCC$ of $\CCC$. Category  $\kk\CCC$ has the same 
objects as $\CCC$, that is, finite sequences $\undep$ of plus and minus signs. Morphisms in $\kk\CCC$ are  finite linear combinations of  morphisms  in  $\CCC$, with the composition rules extended $\kk$-bilinearly from those  of  $\CCC$. 

\vspace{0.1in}  

{\it (2)  Category $\vCCC_{\alpha}$  of viewable cobordisms.}
Next, choose a symmetric power series $\alpha\in\kllS^s$. 
Define the  category $\vCCC_{\alpha}$ as the quotient of 
$\kk\CCC$ by the relations that a circle $\widehat{w)}$ with a sequence $w$ written on it evaluates to $\alpha(w)$. Since $\widehat{w_1w_2}=\widehat{w_2w_1}$, we need the condition that  $\alpha$ is symmetric to define this quotient.  

Another way to define $\vCCC_{\alpha}$ is to say that  it has  the same  objects as $\CCC$: sequences $\undep$ of elements of $S$. A morphism in  $\vCCC_{\alpha}$ from $\undep$  to $\undep'$ is a finite $\kk$-linear combination of viewable cobordisms in $\CCC$ from $\undep$ to $\undep'$. Recall that a cobordism is viewable if it has no floating connected components, that is, components homeomorphic to circles. 

Composition of morphisms in $\vCCC_{\alpha}$ is given by concatenating  cobordisms and removing each  closed circle $\widehat{w}$ from the composition  simultaneously with  multiplying  the remaining diagram by $\alpha(w)$. 

\begin{figure} 
\begin{center}
\includegraphics[scale=1.0]{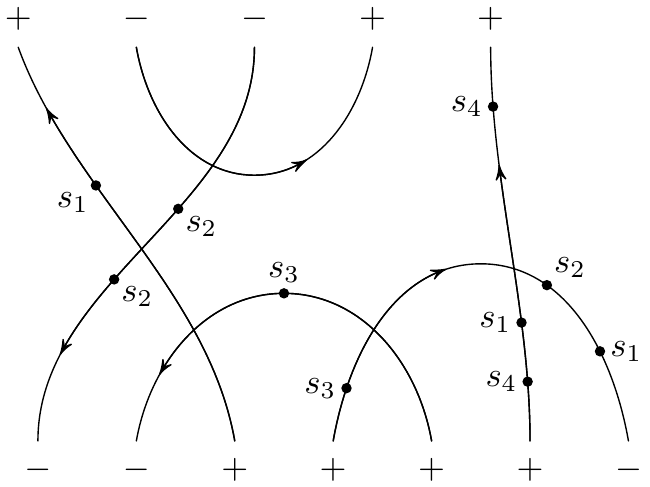}
\caption{A  basis  element in the hom space in $\vCCC_{\alpha}$,  with $S=\{s_1,s_2,s_3,s_4\}$. Floating components (circles) are absent.  }
\label{fig2_1}
\end{center}
\end{figure}

The hom space $\Hom_{\vCCC_{\alpha}}(\undep, \undep')$ has a basis  given by a choice of orientation-respecting matching of the elements in the pair of  sequences $\undep, \undep'$ together with a choice of  word in $S$ for  each pair in the matching. An orientation-respecting  matching consists of a bijection between pluses and minuses in the  sequence $(-\undep)\undep'$, which is the concatenation  of $-\undep$ and $\undep'$, with the sequence $-\undep$  given by reversing the signs  of $\undep$. An example  in  Figure~\ref{fig2_1} shows a basis  element in one such hom space, with $\undep=(--++++-)$ and $\undep'=(+--++)$. Note that the size of hom spaces in $\vCCC_{\alpha}$ does not depend on $\alpha$,  only  the composition of morphisms does. 

\vspace{0.1in}

\begin{figure} 
\begin{center}
\includegraphics[scale=1.0]{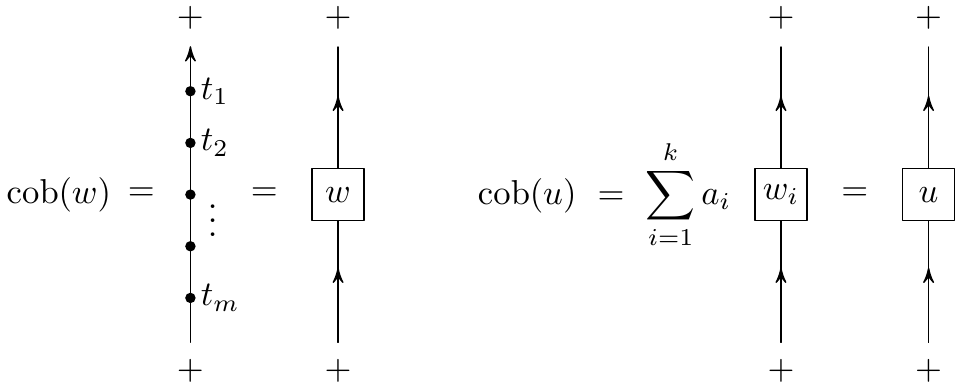}
\caption{ Left: cobordism $\cob(w)$ for a word $w=t_1\dots t_m\in S^{\ast}$  is given by placing  dots labelled by letters of word  $w$ along the oriented interval. Alternatively, $\cob(w)$ can be  denoted by a  box  labelled  $w$ on an interval. Right: a linear  combination $\cob(u)$ of  such cobordisms and  its shorthand box  notation. }
\label{fig2_2}
\end{center}
\end{figure}

Each word $w=t_1\dots t_m,$  $t_i\in S$, $i=1,\dots, m$ defines a cobordism $\cob(w)$ given by putting letters $t_1, \dots, t_m$ along the interval, with  the orientation going towards decreasing the index, see Figure~\ref{fig2_2} left.   Extended by linearity, this assignment  is an algebra isomorphism 
\begin{equation}\label{eq_end}
    \klS \  \lra  \ \End_{\vCCC_{\alpha}}((+))
\end{equation}
from the  algebra of noncommutative polynomials to the  algebra of endomorphisms of the sequence $(+)$ in $\vCCC_{\alpha}$.  

To a noncommutative  polynomial  
\begin{equation*}
u \ = \ \sum_{i=1}^k a_i w_i \in \klS, \ \ a_i\in \kk, \  w_i\in S^{\ast}
\end{equation*} 
we assign the endomorphism 
\begin{equation*}
\cob(u) \ = \ \sum_{i=1}^k a_i \, \cob(w_i) \ \in \  \End_{\vCCC_{\alpha}}((+))
\end{equation*}
of the sequence $(+)$ given  by the linear combination of words $w_i$ written on an upward oriented interval, see Figure~\ref{fig2_2} right. It can  be compactly denoted by  a box on a strand  with $u$ written  in it. 

Monomials in $\klS$ and their linear  combinations can be  placed along any component of a cobordism. Taking the union over all viewable cobordisms with a given  boundary (one cobordism for each diffeomorphism class rel boundary) and then  over all ways of placing  monomials in $\klS$ along each component of  the  cobordism gives a basis in the hom space in the category $\vCCC_{\alpha}$ between two  objects. Recall that objects of $\vCCC_{\alpha}$ are sequences of $+$ and $-$. 

The hom  space between  the objects $\undep$, $\undep'$ is non-zero if the objects  have the same weight, $|\undep|=|\undep'|$. Assuming  the latter, the hom space $\Hom(\undep,\undep')$ is infinite-dimensional unless $\undep=\undep'=\emptyset$ is the empty  sequence or if the set $S$ of labels is empty. In the latter case $\klS\cong \kk$ is the ground field. Another special case is when $S$ consists of a single element, $S=\{s\}$, for then  $\klS\cong  \kk[s]$ is commutative. The endomorphism algebra of $(+)$ in the category $\vCCC_{\alpha}$ is $\klS$, see (\ref{eq_end}).  

\vspace{0.05in} 

Category $\vCCC_{\alpha}$ is  a $\kk$-linear pre-additive category. 

\vspace{0.1in}

{\it (3) The skein  category $\sCCC_{\alpha}$.} 
Consider  the syntactic ideal $I_{\alpha}\subset \klS$ associated to the symmetric series $\alpha\in\kllS^s$. This  ideal has finite  codimension iff $\alpha\in\kllS^{s,rec}$, that is, if $\alpha$ is,  in addition, a recognizable series.
Denote by  
\begin{equation}
    A_{\alpha} \  := \ \klS/I_{\alpha}
\end{equation}
the  quotient algebra (the syntactic algebra) of the algebra of  noncommutative  polynomials by the syntactic ideal. Algebra  $A_{\alpha}$ is finite-dimensional iff $\alpha$ is recognizable.  In the latter case, let 
\begin{equation}
    d_{\alpha} = \dim_{\kk} (A_{\alpha}) =  \mathrm{codim}_{\kk}(I_{\alpha}) 
\end{equation}
be the dimension of the syntactic algebra. 

\vspace{0.1in}

We quotient  the category $\vCCC_{\alpha}$ of viewable cobordisms by the relation that elements of $I_{\alpha}$ are zero along any component  of a cobordism. 
Namely, an element of $I_{\alpha}$ is a finite linear  combination 
\begin{equation}
    u = \sum_{i=1}^k a_i  w_i, \ \ a_i\in \kk, \ w_i\in  S^{\ast}
\end{equation}
of words in the alphabet $S$.  
Element $\cob(u)$, see Figure~\ref{fig2_2}, can be inserted along any component of a cobordism $x$. We  impose the  condition  that any such insertion  results in  the zero morphism in $\sCCC_{\alpha}$ between the corresponding sequences $\undep,\undep'$. Equivalently, we can set $\cob(u)\in  \End((+))$ to zero for all $u\in I_{\alpha}$ and  take the monoidal closure of the relations $\cob(u)=0$ for all such $u$, which is equivalent to the  previous condition. Alternatively, we can choose generators $\{u_j\} , j\in J$, for  the  2-sided  ideal  $I_{\alpha}$, impose relation $\cob(u_j)=0, j\in J$ and  take their monoidal closure. 

\begin{figure} 
\begin{center}
\includegraphics[scale=1.0]{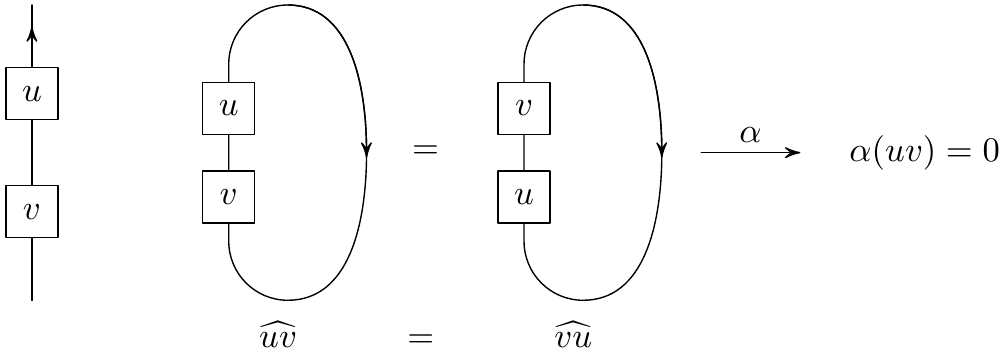}
\caption{$\alpha(\widehat{uv})=\alpha(\widehat{vu})=0$ for $u\in  I_{\alpha}$,  $v\in\klS$  since $uv,vu\in I_{\alpha}$.}
\label{fig2_3}
\end{center}
\end{figure}

Note that relations $\cob(u)=0$ for $u\in I_{\alpha}$ are compatible with the evaluation of closed components (circles). Namely, for any $v\in \klS$, the closures $\widehat{uv}$ and $\widehat{vu}$ define the same element in $\End_{\kk\CCC}(\emptyset)$, namely the circle that  carries the box $uv$ or $vu$, and $\alpha(uv)=\alpha(vu)=0$, see  Figure~\ref{fig2_3}.  Consequently, no contradiction in evaluation of closed components happens upon introducing these relations.

Denote by $\sCCC_{\alpha}$ the resulting quotient category. It has  the same  objects as $\vCCC_{\alpha}$ and additional relations $\cob(u)=0$ for $u\in  I_{\alpha}$ placed  anywhere  along one-dimensional $S$-decorated cobordisms that span hom spaces in $\vCCC_{\alpha}$.

Since relations in the syntactic ideal are imposed along each connected component of a cobordism, an element along a component can be reduced  accordingly. Choose a set of elements $B_{\alpha}\subset \klS$ that descend to a basis of $A_{\alpha}$ (if needed, one  can choose monomials in $S$). Modulo  $I_{\alpha}$, an element of $\klS$ can be  reduced to a  linear  combination of elements of $B_{\alpha}$. Accordingly, we can reduce a morphism in $\sCCC_{\alpha}$ to a linear  combination of viewable morphisms such that along  each component an  element of $B_{\alpha}$ is placed. 
Call these morphisms \emph{basic} and denote the set of  basic morphisms from $\undep$  to $\undep'$ by $B_{\alpha}(\undep,\undep').$ 

Recall that a morphism from $\undep$ to $\undep'$  exists in $\CCC$ if the two sequences have the same weight, that  is, the difference between the number of plus and minus signs in them: $|\undep|=|\undep'|.$ In the latter case, the number of  viewable morphisms (i.e., without circle components) is the  number of ways to pair  up elements of $\undep$ and $\undep'$ in an orientation-respecting way. Reverse the signs in  one  of  the sequences, say in $\undep$,  and concatenate  with the other to get $\undep'(-\undep)$. This sequence  has the same number $n$ of plus and  minus signs, equal to half the length of the sequence: $2n=||\undep||+||\undep'||$. Isomorphism  classes of  viewable  cobordisms from $\undep$ to $\undep'$ are in a one-to-one correspondence with bijections  between  plus and minus signs in $\undep'(-\undep)$. There are $n!$ such bijections. For each bijection, there are $d_{\alpha}^n$ ways to assign an element of $B_{\alpha}$ to each component of a cobordism. The following proposition and corollary result. 

\begin{prop} The set of basic morphisms $B_{\alpha}(\undep,\undep')$ is  a basis  of  the hom  space $\Hom_{\sCCC_{\alpha}}(\undep,\undep')$. 
\end{prop} 

\begin{corollary} Dimensions of hom  spaces in  $\sCCC_{\alpha}$ are given  by: 
\begin{equation}
    \dim \Hom_{\sCCC_{\alpha}}(\undep,\undep')  = \begin{cases} \ \   n!\, d_{\alpha}^{\, n} \  \mbox{ if }  |\undep|=|\undep'|, \ 2n = ||\undep||+||\undep'|| , \\  \ \  0 \quad \quad \  \mbox{otherwise}.
    \end{cases}
\end{equation}
\end{corollary} 
In particular,  hom spaces in the category $\sCCC_{\alpha}$ are finite-dimensional. 

For the endomorphism algebra of the sequence $(+)$ we have (compare with (\ref{eq_end}))
\begin{equation}
    \End_{\sCCC_{\alpha}}((+)) \cong A_{\alpha},
\end{equation}
and the endomorphism algebra of $(+)$ has dimension $d_{\alpha}$. Skein category  $\sCCC_{\alpha}$ is similar to the oriented Brauer category~\cite{R}, but with lines decorated   by elements of $S$, leading to many choices for  evaluations of floating components, one  for each sequence in $S^{\ast}$ up to the cyclic equivalence. 

\vspace{0.1in} 

{\it (4) Negligible morphisms and  gligible quotient $\CCC_{\alpha}$.}

\begin{figure}[h]
\begin{center}
\includegraphics[scale=1.0]{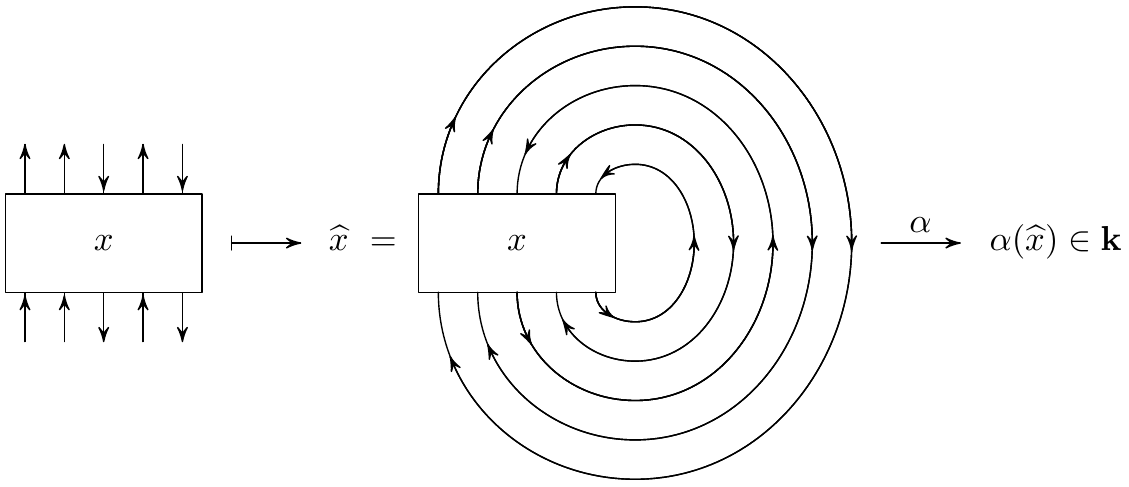}
\caption{The trace map: closing endomorphism $x$ of $\undep$ into  $\widehat{x}$ and applying $\alpha$. In this example $\undep=(++-+-)$.   }
\label{fig3_3}
\end{center}
\end{figure}

The \emph{trace} $\tral(x)$ of a cobordism $x$ from $\undep$ to $\undep$  is an element of  $\kk$  given  by closing $x$ via 
$||\undep||$ suitably oriented arcs  connecting $n$ top with $n$ bottom points of $x$ into a floating cobordism $\widehat{x}$ and applying $\alpha$, 
\begin{equation*}
  \tral(x) := \alpha(\widehat{x}). 
\end{equation*} 
This operation is depicted in Figure~\ref{fig3_3}. The trace is  extended to all endomorphisms  of $\undep$ in  $\kk\CCI$ by linearity. It is  well-defined on trace of endomorphisms of objects  $\undep$ in categories $\vCCI_{\alpha}$ and  $\sCCI_{\alpha}$ as well. 

The trace is symmetric, that is $\tral(yx)=\tral(xy)$ for a morphism $x$ from $\undep$ to  $\undep'$ and $y$ from $\undep'$ to  $\undep$. The ideal $J_{\alpha}\subset \sCCC_{\alpha}$ is defined as follows. 

A morphism $y\in\Hom(\undep,\undep')$ is called \emph{negligible} and belongs to the ideal $J_{\alpha}$ 
if $\tral(zy)=0$ for any morphism $z\in\Hom(\undep',\undep)$.
  Negligible morphisms constitute a two-sided ideal in the pre-additive category $\sCCC_{\alpha}$. 
We call $J_{\alpha}$ the \emph{ideal of negligible morphisms}, relative to the  trace form $\tral$. 
Define the quotient category  
\begin{equation*}
    \CCC_{\alpha} := \sCCC_{\alpha}/J_{\alpha}.
\end{equation*}

 The quotient category $\CCC_{\alpha}$ has finite-dimensional hom spaces, as does  $\sCCC_{\alpha}$ (recall that $\alpha$ is recognizable). The trace is nondegenerate on $\CCC_{\alpha}$ and defines perfect bilinear  pairings
 \begin{equation*}
     \Hom(\undep,\undep') \otimes \Hom(\undep',\undep) \lra \kk
 \end{equation*}
 on its hom spaces. We may call $\CCC_{\alpha}$ the  \emph{gligible quotient} of $\sCCC_{\alpha}$, having modded out by  the ideal of negligible morphisms. 
 
\vspace{0.1in} 

{\it State spaces of recognizable series $\alpha$.}
 Recall that in the category $\vCCC_{\alpha}$ objects are sign sequences $\undep$ and morphisms are finite linear combinations of viewable cobordisms. 
 The space  of homs 
 \begin{equation*}
 V_{\undep} \ := \ \Hom_{\vCCC_{\alpha}}(\emptyset,\undep),
 \end{equation*}
 has a basis of all viewable cobordisms (no floating components) $M$ with $\partial M = \undep$. This space carries a symmetric bilinear  form, given on pairs of  basis elements (viewable cobordisms)  by 
 \begin{equation*}
     (x,y)_{\undep} \ := \alpha(\overline{y}x) \in \kk, 
 \end{equation*}
 where $\overline{y}$ is the reflection of  $y$ about  a horizontal line combined with the orientation reversal on $y$, and $\overline{y}x$ is the closed cobordism which is the  composition of $\overline{y}$ and $x$. 
 
  Define $A_{\alpha}(\undep)$ as the quotient of $V_{\undep}$ by  the  kernel of this  bilinear form. Then there is a  canonical isomorphism
 \begin{equation*}
     A_{\alpha}(\undep) \cong \Hom_{\CCC_{\alpha}}(\emptyset,\undep)
 \end{equation*}
 as  well as isomorphisms 
 \begin{equation*}
     A_{\alpha}((-\undep)\sqcup\undep') \cong \Hom_{\CCC_{\alpha}}(0,(-\undep)\sqcup \undep')  \cong \Hom_{\CCC_{\alpha}}(\undep,\undep')
 \end{equation*}
 given by moving the bottom boundary $\undep$ of a cobordism  to the top  via a cobordism with $||\undep||$ parallel arcs. Here the sequence $(-\undep)\sqcup\undep'$ is the concatenation of $-\undep$ and $\undep'$. 
 
 Note that $\undep$ must be balanced for $A_{\alpha}(\undep)$ to be nonzero, that is, $\undep$ must have the same number $n$ of pluses and minuses. We can  then define
 \begin{equation}\label{eq_A_n} 
     A_{\alpha}(n) \ := \ A_{\alpha}((+^n-^n)).
 \end{equation}

 Spaces $A_{\alpha}(n)$ come with a lot of structure, including  multiplication  maps 
 \begin{equation*} 
 A_{\alpha}(n) \otimes  A_{\alpha}(m)\lra A_{\alpha}(n+m). 
 \end{equation*} 
 
 We have $A_{\alpha}(0)\cong \kk$ and $A_{\alpha}(1)\cong A_{\alpha}$.  Vector space $A_{\alpha}(n)$ carries an  action  of the  symmetric  group product $S_n\times S_n$ by the permutation cobordisms, as well as an action of the tensor power of the syntactic algebra $A_{\alpha}^{\otimes n} \otimes (A_{\alpha}^{op})^{\otimes n}$, with one copy of $A_{\alpha}\cong \End_{\CCC_{\alpha}}((+))$ or $A_{\alpha}^{op}\cong \End_{\CCC_{\alpha}}((-))$   acting at each sign of $+^n-^n$. More generally, a version  of the oriented walled Brauer algebra with strands carrying $S$-labelled dots and closed decorated circles evaluating via  $\alpha$ acts on $A_{\alpha}(n)$ and, more generally, on $\Hom_{\CCC_{\alpha}}(\undep,+^n-^n)$ for any sign sequence $\undep$. This generalized walled Brauer algebra $\mathrm{Br}_{n,\alpha}$ is straightforward to define. It is associated to any recognizable series $\alpha$, finite-dimensional, and isomorphic to the endomorphism algebra 
 $\End_{\sCCC_{\alpha}}((+^n-^n))$ of the object $(+^n-^n)$ in the skein category $\sCCC_{\alpha}$. The action of $\mathrm{Br}_{n,\alpha}$ on $\Hom_{\CCC_{\alpha}}(\undep,(+^n-^n))$ descends to the action of its quotient algebra    $\End_{\CCC_{\alpha}}((+^n-^n))$ on the same space.
 
 Multiplication maps turn the direct  sum 
 \begin{equation}\label{eq_A_n_1} 
     A_{\alpha}^{\ast} \ := \ \oplusop{n\ge 0} A_{\alpha}(n)
 \end{equation}
 into  a graded associative $\kk$-algebra, with compatible actions of $S_n$ on $A_{\alpha}(n)$ over all $n$, making $A_{\alpha}^{\ast}$ into what Sam and Snowden call a \emph{twisted commutative algebra} or \emph{tca} in~\cite[Definition 7.2.1]{SSn}. A twisted commutative algebra in that sense may be very from being  commutative: for instance, the free  associative  algebra (the tensor  algebra of a vector space) has the  obvious tca structure~\cite[Example 7.2.2]{SSn}.
 More generally, given an $n$-dimensional  topological theory $\alpha$ as  defined in~\cite{Kh2}, perhaps  for  manifolds with defects, etc. and an $(n-1)$-manifold $N$, the direct sum 
 \begin{equation*}
     A^{\ast}(N) \ := \ \oplusop{n\ge 0} \alpha(\sqcup_n N) 
 \end{equation*}
 of state spaces  of disjoint unions of $n$ copies of $N$, over all $n$, is naturally a tca in the sense of~\cite{SSn}. 
 
\vspace{0.1in} 

{\it (5) The Deligne category $\dCCC_{\alpha}$ and its gligible  quotient $\udCCC_{\,\alpha}$.}
The skein category  $\sCCC_{\alpha}$  is a rigid symmetric monoidal $\kk$-linear category with signed sequences $\undep$ as objects and finite-dimensional hom spaces. 
We form the additive Karoubi closure 
\begin{equation*}
    \dCCC_{\alpha} \ := \ \Kar(\sCCC^{\oplus}_{\alpha})
\end{equation*}
by allowing formal finite direct sums of objects in $\sCCC$,   extending morphisms correspondingly, and then adding idempotents to get a Karoubi-closed  category. Category $\dCCC_{\alpha}$ plays the role of the Deligne category in  our construction.  

\vspace{0.1in} 

The trace $\tral$ extends to $\dCCC_{\alpha}$ and defines a 2-sided ideal $\mathcal{DJ}_{\alpha}\subset \dCCC_{\alpha}$ of negligible morphisms relative to $\tral.$ Define the \emph{gligible quotient} category by 
\begin{equation*}
    \udCCC_{\alpha} \ := \ \udCCC_{\alpha}/ \mathcal{DJ}_{\alpha}.
\end{equation*}
This category is equivalent to the additive Karoubi envelope of $\CCC_{\alpha}$. It is a  Karoubi-closed rigid symmetric category with non-degenerate bilinear forms on its hom spaces. 
\vspace{0.1in} 


\subsection{Summary of categories  and  functors}\label{subsec_summ}
$\quad$ 
\vspace{0.1in} 

Here is the summary of the categories  that have been  introduced.  

\begin{itemize}
    \item  $\CCC$: the category of $S$-decorated one-dimensional  cobordisms. Its objects are sequences $\undep$ of plus and minus signs and morphisms are one-manifolds  with boundary decorated by $S$-labelled dots. That is, the  morphisms are $S$-decorated one-manifolds with boundary. 
    \item $\kk\CCC$: this category has the  same objects as $\CCC$; its  morphisms are formal finite $\kk$-linear combinations of morphisms in $\CCC$. 
    \item  $\vCCC_{\alpha}$: in this  quotient category of $\kk\CCC$ we reduce morphisms to linear combinations  of  viewable cobordisms. Floating  connected components (circles, possibly carrying $S$-dots) are  removed  by  evaluating them  via  $\alpha$.  
    \item $\sCCC_{\alpha}$: to define this category,   specialize to rational $\alpha$ and add skein relations by modding out  by elements of the ideal  $I_{\alpha}$ in  $\kk[S]$, along each component of the cobordism. Hom spaces in this category are finite-dimensional. 
    \item  $\CCC_{\alpha}$: the quotient  of  $\sCCC_{\alpha}$ by the  ideal $J_{\alpha}$ of  negligible morphisms.  This category  is also equivalent  (even isomorphic) to the quotients   of $\kk\CCC$ and $\vCCC_{\alpha}$ by the corresponding ideals of negligible morphisms in  them. The trace pairing in $\CCC_{\alpha}$ between $\Hom(n,m)$ and $\Hom(m,n)$ is perfect. 
    \item $\dCCC_{\alpha}$ is the  analogue  of the Deligne category obtained from $\sCCC_{\alpha}$ by allowing finite direct sums of objects and then  adding  idempotents  as objects to get a Karoubi-closed category. 
    \item $\udCCC_{\,\alpha}$: the quotient  of $\dCCC_{\alpha}$ by the two-sided ideal of negligible  morphisms. This category is equivalent to the additive  Karoubi closure of $\CCC_{\alpha}$ and sits in  the bottom right corner of the  commutative square below.  
\end{itemize}

We arrange these  categories and functors, for   recognizable $\alpha$,
into  the following diagram: 

\begin{equation} \label{eq_seq_cd_1}
\begin{CD}
\CCC  @>>> \kk\CCC @>>> \vCCC_{\alpha} @>>> \sCCC_{\alpha} @>>>  \dCCC_{\alpha} \\
@.   @. @.   @VVV   @VVV  \\    
 @.  @.   @.  \CCC_{\alpha}  @>>>  \udCCC_{\,\alpha}
\end{CD}
\end{equation}

All seven  categories  are rigid symmetric monoidal. All but the leftmost  category $\CCC$ are  $\kk$-linear. Except for the two categories on the far right, the  objects of  each category are sequences $\undep$  of plus and minus signs. 
The four categories on the right  all  have finite-dimensional hom spaces. The  two categories on the far right are additive and Karoubi-closed. The four categories in the middle of the  diagram  are pre-additive but not additive.  

The arrows show functors between these   categories considered  in the paper. The square is commutative. An analogous diagram of functors and categories can be found in~\cite{KS3} for  the  category of oriented 2D cobordisms in  place of $\CCC$ and in~\cite{KQR} for suitable  categories of oriented 2D cobordisms with side boundary and corners. 

For convenience, one- or two-word summaries of these categories are provided below, in the diagram essentially identical to that in~\cite[Section 3.4]{KQR}: 

\begin{equation}\label{eq_words}
\begin{gathered}
 \xymatrix{
\vcenter{\hbox{$S$\textrm{-dotted}}\hbox{\textrm{cobordisms}}}
\ar[r] & \kk\textrm{-linear} \ar[r] & \textrm{viewable} \ar[r] & \textrm{skein}\ar[d] \ar[r] &
{\textrm{Deligne (Karoubian)}} \ar[d] \\
& & & \textrm{gligible} \ar[r] & 
\textrm{gligible and Karoubian}
}
\end{gathered}
\end{equation}

It is possible to go directly from $\kk\CCC$ to  $\CCC_{\alpha}$ by modding  out by  the ideal of negligible morphisms  in the former category. It is convenient to  arrive at this  quotient  in several steps, introducing  categories  $\vCCC_{\alpha}$ and  $\sCCC_{\alpha}$ on the way. 

\vspace{0.1in}

If $\alpha$  is not recognizable, we can still  define categories $\vCCC_{\alpha}$, $\CCC_{\alpha}$ and $\udCCC_{\,\alpha}$, but then, for instance,  one  can potentially get two non-equivalent categories in  place of $\udCCC_{\alpha}$ by  following along the two different paths in the square  above. To justify  considering these categories for some non-recognizable $\alpha$ one would want  to find interesting examples where the gligible quotient category $\CCC_{\alpha}$ has additional relations beyond 
those in $\sCCC_{\alpha}$, that is, beyond the relations that elements of the  syntactic ideal $I_{\alpha}$ are zero in $\End(+)$ in $\sCCC_{\alpha}$ and $\CCC_{\alpha}$. 


\subsection{Examples and variations of the construction} 
\quad
\vspace{0.1in} 

{\it An involution.} Categories $\CCC$ and $\kk\CCC$ carry   contravariant  involution $\ophana$ that reflects the cobordism about the middle, reversing its source and target objects, and reverses the orientation of the cobordism. 
This involution takes the object $\undep$ to $-\undep$, that is, reverses the sign (orientation) of boundary zero-manifolds as well. To  match this  involution to evaluation $\alpha$, assume that $\kk$ comes with an involution, also denoted 
$\ophana$, and $\alpha$ satisfies $\overline{\alpha(w)}=\alpha(\overline{w})$, where $\overline{w}=t_n\dots t_1$ is  the word $w=t_1\dots t_n$ in reverse. 
Then there  are induced contravariant involutions on all the  categories associated to $\alpha$ and displayed  in diagram (\ref{eq_seq_cd_1}), and  one can, for  instance,  study such unitary 1D topological theories, with the set of defects $S$, for  $\kk=\C$ and $\ophana$ the complex involution. 

\vspace{0.1in} 

\vspace{0.1in} 

{\it Examples.}

(1)
If the set  $S=\emptyset$ is empty, there are no decorations and the series $\alpha$ is given  by its  value on the empty sequence, that is, by its  constant term, and we can  write $\alpha=\lambda \in  \kk$  for that  value. A circle cobordism evaluates  to $\lambda$. The skein category  $\sCCC_{\lambda}$ is isomorphic to the viewable category $\vCCC_{\lambda}$  and to the oriented  Brauer category  $B_{\lambda}$  for  for the parameter $\lambda$. Category $\CCC_{\lambda}$ is then  the  quotient of $B_{\lambda}$ by the ideal of  negligible morphisms, while $\dCCC_{\lambda}$ is the additive Karoubi closure of $B_{\lambda}$, etc. Note that these categories depend both on the field $\kk$  and  $\lambda\in \kk$. 

\vspace{0.1in}

(2) 
If $S=\{s\}$ is a one-element set, the series $\alpha$ is a one-variable series, with   the generating  function 
\begin{equation}
Z_{\alpha}(T)\ =\ \sum_{n\ge 0} \alpha_n T^n, \ \ \alpha_n=\alpha(s^n).  
\end{equation} 
$\alpha$ is recognizable iff $Z_{\alpha}(T)$ is a rational function, with $Z_{\alpha}(T)=P(T)/Q(T)$ for some polynomials $P(T),Q(T)$. 

This example is similar to the ones in~\cite{Kh2,KS3}, where the topological theory is 2-dimensional but there are no defects. The analogue of the Hankel matrix measuring bilinear pairing on connected cobordisms with the boundary $\SS^1$ (such cobordisms are determined by the genus $g$), see~\cite{Kh2}, is the Hankel matrix for evaluations of $\overline{x_m}x_n$ where $x_n$ is an arc with $n$ dots, viewed as a 
cobordism from $\emptyset$ to $(+-)$. Cobordism $\overline{x_m}$ from $(+-)$ to $\emptyset$ is an arc with $m$ dots. Closed cobordism  $\overline{x_m}x_n$ is a circle with $n+m$ dots, and once again the Hankel matrix $H$ with the $(n,m)$-entry $\alpha_{n+m}$ results, as in~\cite{Kh2}. In both cases the state space (of $(+-)$, respectively of $\SS^1$) is the quotient of  $\R^{\N}$ by the null space of $H$. 

The theories diverge beyond this  example, but there is another connection between  the two, slightly  different from the one above due to an additional shift in the dots versus handles correspondence between 1D and 2D cobordisms. Namely, the  state space of $(+^k-^k)$ in the one-dimensional theory with the series  $\alpha$ maps to  the state space for the union $\sqcup_k \SS^1$ of $k$ circles in the  two-dimensional theory for the series $\alpha'=(a,\alpha_0,\alpha_1,\dots )$ for any $a\in \kk$. In terms of generating functions, $Z_{\alpha'} = a + Z_{\alpha} T$. On  the topological side, an arc with $n$ dots is mapped to a  an annulus with $n$ handles, while a circle carrying $n$ dots is mapped to to the torus with additional $n$ handles (thus a  surface of genus $n+1$). This shift from $n$ to $n+1$ accounts for the discrepancy between the series but does not change their recognizability. The  map from the state spaces in the 1D theory to the state spaces in the 2D theory  respects the bilinear forms on these spaces. 

Partial fraction decomposition method of~\cite{KKO} can be applied in this case as well to understand the categories associated to $\alpha$. When the set $S$ has more than  one element, recognizable power series still admit an analogue of the partial  fraction decomposition, see~\cite{FH} and references  therein, which should lead to decompositions of associated tensor categories. 

\vspace{0.1in} 


{\it Unoriented cobordisms.} There is an  obvious unoriented  version of the category  $\CCC$, where one-dimensional cobordisms are  unoriented and the objects, in the  skeletal category case, are numbers  $n\in \Z_+$, counting  the number of  top and  bottom endpoints of the  cobordism. Evaluation  $\alpha$ must  be  $\ophana$-invariant, that is,  to satisfy $\alpha(\overline{w})=\alpha(w)$,  for any  word  $w$, in addition to being  symmetric, as earlier: $\alpha(w_1w_2)=\alpha(w_2w_1)$, for any words $w_1,w_2$. The dihedral group  $D_n$  acts  on the set $S^n$ of  words of length  $n$ in the alphabet  $S$, and the  function $\alpha:S^{\ast}\lra \kk$, when restricted to these words, must  be $D_n$-invariant. Such series can be called  \emph{d-symmetric}, for instance. 

The theory  then goes  through  and  one can define the viewable category $\vCCC_{\alpha}$, the skein category $\sCCC_{\alpha}$, the  gligible quotient $\CCC_{\alpha}$, and so  on. The interesting case, as before, is when $\alpha$ is \emph{recognizable}, that is, when the category $\CCC_{\alpha}$ has finite-dimensional hom spaces. A d-symmetric series  $\alpha$ is recognizable iff it is recognizable as noncommutative series iff  the syntactic  ideal $I_{\alpha}$ has finite codimension in $\klS$. 

If the set $S$ is empty, cobordisms do not carry any dots (defects), and  the category $\sCCC_{\alpha}$ is the unoriented  Brauer category $\mathrm{Br}^{un}_{\lambda}$ for the parameter $\lambda=\alpha(1)\in\kk$, while $\CCC_{\alpha}$ is the gligible quotient of $B^{un}_{\lambda}$. 

\vspace{0.1in} 

%
%

\section{Cobordisms with inner (floating) boundary }
\label{sec_floating_cobs}

%
%

\subsection{Category  \texorpdfstring{$\CCI$}{C} of decorated  cobordisms with inner  boundary}
\quad
\vspace{0.1in} 

To connect decorated one-dimensional cobordisms with noncommutative rational power series that are not necessarily symmetric we  enlarge the category $\CCC$ by allowing cobordisms $M$ that may have additional boundary points (\emph{floating boundary} points) strictly  inside the cobordism,  not being part of the top $\partial_1 M$ or bottom $\partial_0 M$  boundary of $N$. 

\begin{figure} 
\begin{center}
\includegraphics[scale=1.0]{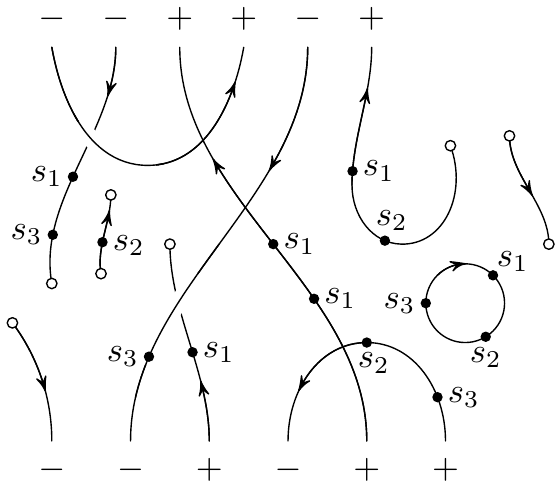}
\caption{A morphism from $(--+-++)$ to $(--++-+)$. Over- and  undercrossing  and  intersections are "virtual" and should be ignored. Hollow dots are not labels and  show inner (floating) endpoints of the cobordism.}
\label{fig4_1}
\end{center}
\end{figure}

Define the category $\CCI$ of  $S$-labelled  cobordism with floating (or inner) boundary to  have the same objects as $\CCC$, that is, finite sequences  $\undep$  of  plus and  minus signs. A morphism in $\CCI$ from $\undep$ to  $\undep'$, see Figure~\ref{fig4_1} for an example, is a compact oriented $S$-decorated one-manifold $M$ with 
\begin{equation}\label{eq_partial_inner}
    \partial M \ = \ \undep' \sqcup (- \undep)  \sqcup  \partial_{in} M,
\end{equation}
where $\partial_{in} M$ is the \emph{inner} or \emph{floating} boundary  of $M$ that is disjoint from the top boundary, given  by $\undep'$ and from the bottom boundary, given by $-\undep$. In (\ref{eq_partial_inner}) we interpret a sign sequence $\undep'$ as a zero-dimensional oriented manifold, with oriented connected components described by elements of the sequence. The sequence  $-\undep$ opposite to $\undep$ corresponds to the orientation reversal of  0D manifold $\undep$. An  $S$-decoration is a  collection of points (dots) labelled by elements of the set $S$ inside $M$ (not on the  boundary $\partial M$).  Labelled points can move along a connected component but not cross through each other. 

Morphisms are such  decorated 1D cobordisms, possibly with inner endpoints (inner boundary  points)  considered up to rel boundary  diffeomorphisms. Figure~\ref{fig4_1}  shows an example of a morphism from $\undep=(--+-++)$ to $\undep'=(--++-+)$. 

Composition of morphisms  in $\CCI$ is given by concatenation of cobordisms. The category $\CCI$ contains $\CCC$ as the subcategory with  the same objects as $\CCI$ and morphisms -- morphisms of $\CCI$ that have no inner (floating) boundary points. 

Connected components of a cobordism in $\CCI$ split into \emph{viewable} and \emph{floating} types. Figure~\ref{fig4_1} cobordism has three floating components: one circle and  two intervals. The same cobordism has eight viewable components: four  of them have both endpoints on top or bottom boundary, while the other four have one floating endpoint.
Floating components terminology  was introduced in~\cite{KS1}. 

\vspace{0.1in} 

Going along a component $c$ in the  direction of its orientation we read off the labels of dots. If the component is an arc, the result is a sequence $\mathrm{sec}(c)\in S^{\ast}$, a word in the alphabet $S$. If the  component is a circle, the  sequence $\mathrm{sec}(c)$ is defined up to cyclic rotation. Our convention is to write the sequence from right to left as we follow the orientation. For instance, in Figure~\ref{fig4_6} left the sequence is $s_1s_2s_3$, while in Figure~\ref{fig4_6} right the sequence is $s_3s_2s_1$. Orientation reversal of a component corresponds to reversing the sequence. 

\begin{figure} 
\begin{center}
\includegraphics[scale=1.0]{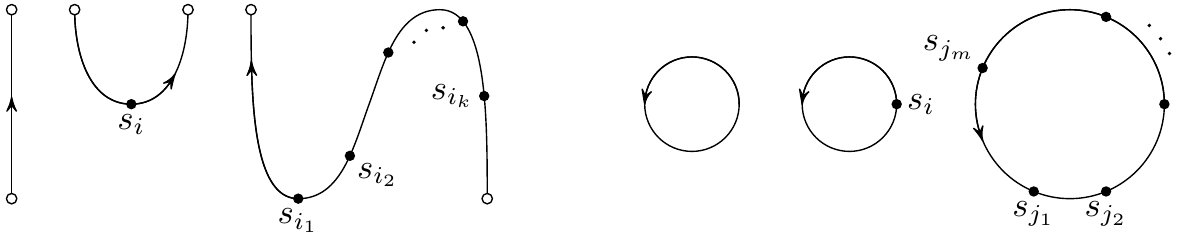}
\caption{Left: three interval floating components, with sequences $(\emptyset)$, $(s_i)$, and $(s_{i_1}s_{i_2}\dots s_{i_k})$. Right: three circle components, with sequences  $(\emptyset)$, $(s_i)$ and $(s_{j_m}\dots s_{j_2}s_{j_1})$ up to cyclic rotation.}
\label{fig4_2}
\end{center}
\end{figure}

The sequences for components of  Figure~\ref{fig4_1} cobordism are: 
\begin{itemize}
    \item The  empty sequence $(\emptyset)$ and $(s_2)$ for the two floating arc components.
    \item Sequence $(s_1s_3s_2)$, up to cyclic rotation, for the unique  floating circle component.
    \item Sequences $(s_1s_1) $ and $(s_3)$ for the two
    connected components that connect a  top endpoint and a bottom  endpoint. 
    \item The empty sequence $(\emptyset)$ for the unique component that connects two top endpoints. 
    \item Sequence $(s_2s_3)$ for the  unique component connecting two bottom endpoints. 
    \item Sequences $(s_3s_1)$ and $(s_1s_2)$ for arc components with one top and one inner boundary point. 
    \item Sequences $(\emptyset)$ and $(s_1)$ for arc components with one bottom  and one  inner boundary  point. 
\end{itemize}
A floating  component of a cobordism $x$ in $\CCI$ is either an interval or a circle, see Figure~\ref{fig4_2}.  
 
\begin{figure} 
\begin{center}
\includegraphics[scale=1.0]{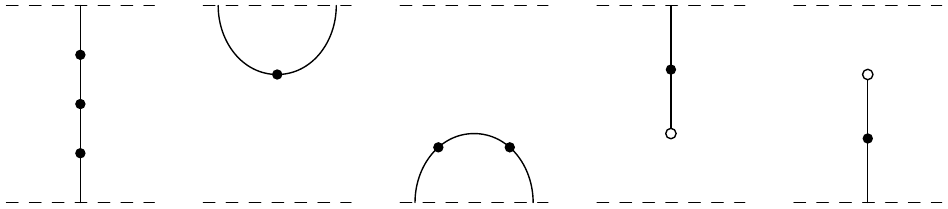}
\caption{Five types of viewable components, left to right: an interval connecting (1) a top and a bottom point, (2) two top  points, (3) two bottom points; a interval with an inner boundary point and a (4) top endpoint, (5) bottom endpoint. Labels of dots and orientations of lines are not shown. In the subcategory  $\CCC$ viewable components are of types (1)-(3) only.}
\label{fig4_3}
\end{center}
\end{figure}
A viewable component has one of the five types shown in Figure~\ref{fig4_3}, with some number of dots (perhaps none) on it.

\begin{figure} 
\begin{center}
\includegraphics[scale=1.0]{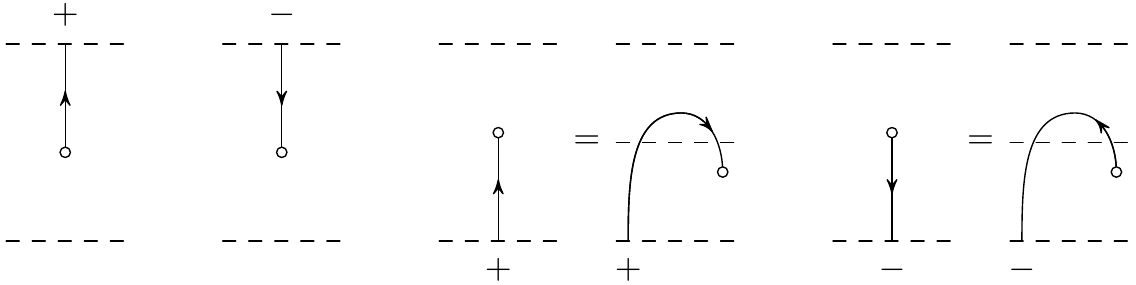}
\caption{Left: additional generating morphisms for monoidal category $\CCI$ beyond the generating morphisms common for $\CCI$ and $\CCC$ shown in Figures~\ref{fig1_3} and~\ref{fig1_4}. }
\label{fig4_4}
\end{center}
\end{figure}

\begin{figure} 
\begin{center}
\includegraphics[scale=1.0]{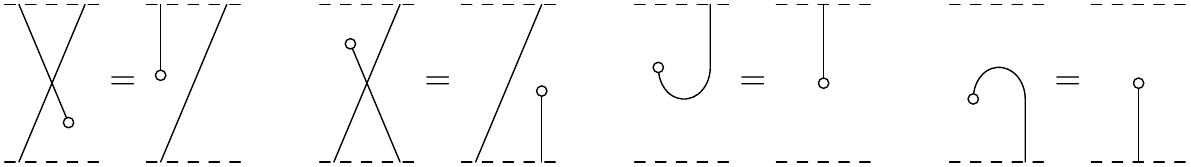}
\caption{Some additional relations in $\CCI$. }
\label{fig4_5}
\end{center}
\end{figure}

\begin{figure} 
\begin{center}
\includegraphics[scale=1.0]{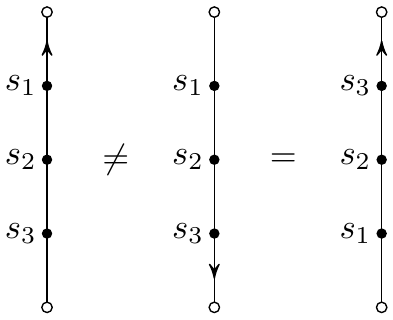}
\caption{Orientation matters: evaluations $\alpha(s_1s_2s_3)$ and $\alpha(s_3s_2s_1)$ are different, in general. Reversal of a sequence corresponds to orientation reversal of the  corresponding floating arc or circle.}
\label{fig4_6}
\end{center}
\end{figure}

Monoidal category $\CCI$ has generators shown in  Figures~\ref{fig1_3} and~\ref{fig1_4} and common with its subcategory $\CCC$ and two additional 
generators shows in Figure~\ref{fig4_4} left. These are ars with one floating  and  one top endpoint. Applying U-turns to them results in arcs with one floating and one bottom endpoint, see Figure~\ref{fig4_4} right. 
Some additional defining relations in $\CCI$ are shown in Figure~\ref{fig4_5}, see also Figure~\ref{fig1_4} for defining relations in the subcategory  $\CCC$, which also give a subset of defining relations in $\CCI$. We will  not need a full set of defining relations for $\CCI$ in this paper. 

\subsection{Tensor envelopes of \texorpdfstring{$\CCI$}{C}}
\quad
\vspace{0.1in} 
 
 Floating (closed) cobordisms in $\CCI$ (endomorphisms of the empty zero-manifold  $(\emptyset)$) are unions of floating intervals and circles. A  floating interval carries  a sequence $w\in S^{\ast}$, a floating circle carries a sequence $v$ well-defined up  to cyclic rotation, $v_1v_2\equiv v_2 v_1$. Consequently a  multiplicative  evaluation of floating  cobordisms in $\CCI$, as  explained in the introduction,  consists  of a pair of series 
 \begin{equation} \label{eq_pair_ser}
 \alpha \ = \  (\alphai,\alphac),
 \end{equation} 
 where $\alphai\in \kllS$ is  a noncommutative series and $\alphac\in \kllS^s$ is a symmetric series. 
 
\begin{figure} 
\begin{center}
\includegraphics[scale=1.0]{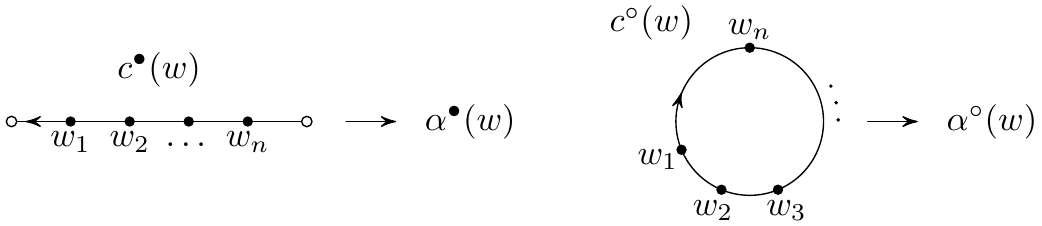}
\caption{Evaluation $\alphai(w)$ of the floating interval $c^{\bullet}(w)$ and evaluation $\alphac(w)$ of the circle $c^{\circ}(w)$ in $\vCCI$, for a word $w=w_1\dots w_n$.}
\label{fig4_7}
\end{center}
\end{figure}

 A multiplicative evaluation on  closed  cobordisms in $\CCI$ assigns    $\alphai(w)\in\kk$ to an oriented  interval with word $w$ written along it via labelled dots, see  Figure~\ref{fig4_7}. Element $\alphac(v)\in \kk$ is assigned to an oriented circle with word $v$, well-defined up to a cyclic rotation, written along it. 

\vspace{0.1in} 

We  now proceed along a  familiar  route, as in~\cite{KS3,KQR} and Section~\ref{sec_decorated}, to 
build  various tensor  envelopes  of  a pair $\alpha=\alphap$. 

\vspace{0.1in} 

(1) {\it  Pre-linearization category  $\kk\CCI$.} Category $\kk\CCI$ has the same  objects as $\CCI$, and the morphisms are finite $\kk$-linear  combinations of morphisms in $\CCI$. This is a naive  \emph{linearization} or  $\emph{pre-linearization}$ of  $\CCI$. 

\vspace{0.1in} 

(2) {\it Viewable cobordisms category $\vCCI_{\alpha}$.} 
To form category $\vCCI_{\alpha}$, we mod out tensor category $\kk\CCI$ by relations that evaluate  floating (closed) cobordisms to elements of the  ground field via $\alpha$. Namely, a floating oriented interval with a sequence $w\in S^{\ast}$ on it,  denoted  $c^{\bullet}(w)$, evaluates to $\alphai(w)\in\kk$. A floating oriented $w$-decorated circle $c^{\circ}(w)$ evaluates to  $\alphac(w)$, see Figure~\ref{fig4_7}.  Recall that $\alphac$ is symmetric and $\alphac(v_1v_2)=\alphac(v_2v_1)$ for  any words  $v_1v_2$, matching circle rotation, $c^{\circ}(v_1v_2)=c^{\circ}(v_2v_1).$

Since all floating  components of a cobordism reduce to elements  in $\kk$, the vector space of homs from  $\undep$ to $\undep'$ in $\vCCI_{\alpha}$  has a basis  of viewable cobordisms from $\undep$ to  $\undep'$ with any sequences written  on  its connected components. 

\vspace{0.1in} 

Denote by $\mcI(\undep,\undep')$ the set of diffeomorphism classes of viewable cobordisms (without dot decorations) from $\undep$ to $\undep'$. 
A viewable cobordism has no circles  and  all its connected components are  intervals.
Such a cobordism $C$ may have some number of viewable  components of types (4) and  (5), see Figure~\ref{fig4_3}. Each such component has one floating boundary  point and  one boundary point among elements of $\undep\sqcup \undep'$. Other  connected components (of types  (1)-(3)) give an orientation-respecting  matching of the remaining elements  of $\undep$ and $\undep'$. 

\begin{figure} 
\begin{center}
\includegraphics[scale=1.0]{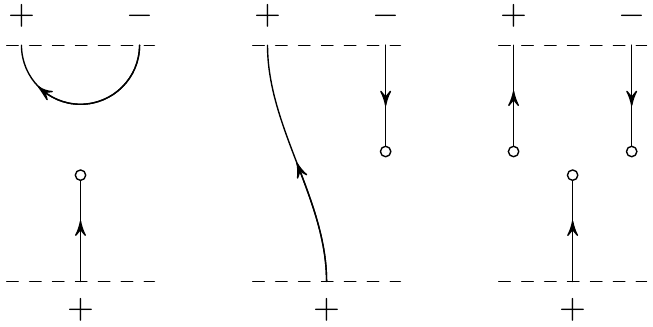}
\caption{Three elements of the set $\mcI(+,+-)$.}
\label{fig4_8}
\end{center}
\end{figure}

\begin{figure} 
\begin{center}
\includegraphics[scale=1.0]{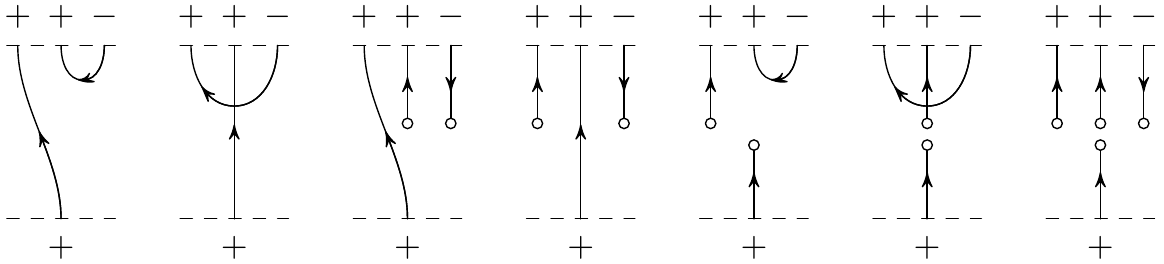}
\caption{Seven elements  of  the set $\mcI(+,++-)$}
\label{fig4_9}
\end{center}
\end{figure}

  To specify an element of $\mcI(\undep,\undep')$ we select a subset $I'$ of elements in the  sequence $(-\undep)\sqcup \undep'$ so that the remaining sequence is \emph{balanced}, that  is, has  the same number of pluses and minuses. We then choose a bijection $b$ between pluses and  minuses  of $(-\undep)\sqcup \undep\, \setminus I'$. Such  pairs  $(I',b)$ are in a bijection with isomorphism classes of  viewable undecorated cobordisms between $\undep$ and  $\undep'$, that  is, elements of $\mcI(\undep,\undep')$. Figure~\ref{fig4_8} shows   elements of  the set $\mcI(+,+-)$.  Figure~\ref{fig4_9} shows   elements of  the set $\mcI(+,++-)$. 
  
  To allow $S$-decorations, we consider the set $\mcI^S(\undep,\undep')$ which consists of a pair: an  element of $\mcI(\undep,\undep')$ and  a  choice  of word $w(c)$ in  $S^{\ast}$ for each component  $c$ of $\mcI(\undep,\undep')$. To such a pair we assign an $S$-decorated viewable cobordism given by the element  of  $\mcI(\undep,\undep')$ and words $w(c)$ written on components $c$ of the cobordism. Some words may be empty (have length zero). Denote the cobordism associated  with $t\in \mcI^S(\undep,\undep') $ by $C(t)$. 

\begin{prop}  Viewable cobordisms $C(t)$, over all $t$ in $\mcI^S(\undep,\undep')$, constitute a basis in the hom space $\Hom(\undep,\undep')$ in the category $\vCCI_{\alpha}$. 
\end{prop}

Composing cobordisms from these bases sets results in cobordisms that, in general, have floating components. These components are evaluated via  $\alpha$, viewable components are kept, and the composition of two basis elements is a basis element in a suitable hom space, scaled by an element of $\kk$. 

\vspace{0.1in}

As earlier, to each word $w\in S^{\ast}$ we associate the upward interval with $w$ written  on  it, that is, cobordism $\cob(w)$ from $(+)$ to $(+)$, see Figure~\ref{fig2_2} left. 
Each element $u$ of $\klS$ gives rise to a linear combination $\cob(u)$ of these cobordisms, see Figure~\ref{fig2_2} right. The  resulting map  
\begin{equation*}
    \cob \ : \ \klS \lra \End_{\vCCI_{\alpha}}((+))
\end{equation*}
is  an injective  homomorphism from the free algebra $\klS$ to the ring of endomorphisms of $(+)$ in category $\vCCI_{\alpha}$ (in the smaller category $\vCCC_{\alpha}$ considered in Section~\ref{sec_decorated} this map is an isomorphism). One-sided inverse homomorphism to $\cob$ is given by the surjection 
\begin{equation*}
    \cob' \ : \  \End_{\vCCI_{\alpha}}((+))\lra \klS 
\end{equation*}
that sends any  cobordism with  floating endpoints to zero. The latter cobordisms span  a two-sided ideal in $\End_{\vCCI_{\alpha}}((+))$, with the quotient isomorphic to $\klS$. This ideal is naturally  isomorphic to $\klS\otimes \klS$ when viewed as a $\klS$-bimodule. Multiplication in this  ideal is given by 
\begin{equation*}
    (x_1\otimes x_2)(y_1\otimes y_2) = \alphai(x_2y_1) \, x_1\otimes y_2. 
\end{equation*}
Composition $\cob'\circ \cob = \mathrm{Id}_{\klS}$. 

\vspace{0.1in} 

(3) {\it Skein category $\sCCI_{\alpha}$.} This category has finite-dimensional hom spaces when  $\alpha$ is recognizable, and we restrict to that case. We say that  $\alpha=\alphap$ is \emph{recognizable} if both series $\alphai$ and $\alphac$ are recognizable. Series $\alphai$ and $\alphac$ has syntactic ideals $I_{\alphai}, I_{\alphac}\subset \klS$, respectively. Recognizability means that both ideals have  finite codimension in $\klS$. Equivalently, the two-sided ideal 
\begin{equation}
    I_{\alpha}  \ := \ I_{\alphai} \cap I_{\alphac} \ \subset \  \klS
\end{equation}
has finite codimension in $\klS$. Denote by 
\begin{equation*}
    A_{\alpha}  \ := \ \klS/I_{\alpha}
\end{equation*}
the syntactic algebra of the pair  $\alpha$. 

\begin{figure} 
\begin{center}
\includegraphics[scale=1.0]{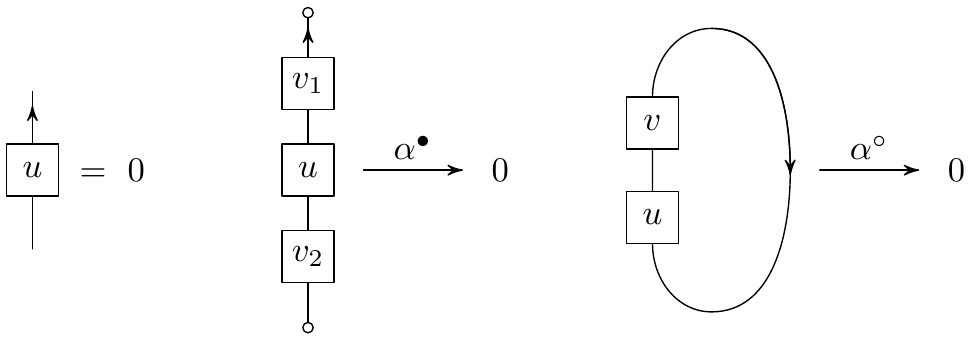}
\caption{Left: endomorphism $\cob(u)$ of $(+)$ is  set to zero  in $\sCCI_{\alpha}$ for $u\in I_{\alpha}$. Middle: $\alphai(v_1uv_2)=0$ for any $v_1,v_2\in\klS$ since $u\in I_{\alpha}\subset I_{\alphai}$. Right: $\alphac(vu)=0$ for any $v\in \klS$ since $u\in I_{\alpha}\subset I_{\alphac}$. }
\label{fig4_10}
\end{center}
\end{figure}

Starting with the category $\vCCI_{\alpha}$, we add tensor relations $\cob(u)=0$ for any $u\in I_{\alpha}$, see Figure~\ref{fig4_10} left. These relations are consistent with the evaluation $\alpha$ of floating components. Consistency is due to restricting to $u$ in the syntactic ideal, which  is contained in both  ideals $I_{\alphai}$ and $I_{\alphac}$. Elements of the first ideal evaluate to zero when placed anywhere on a floating interval, see Figure~\ref{fig4_10} middle. Elements of the second ideal evaluate to zero when placed on a circle, see Figure~\ref{fig4_10} right.  

Recall that in addition to two-sided syntactic ideals $I_{\alphai},I_{\alphac}$ and their intersection $I_{\alpha}=I_{\alphai}\cap I_{\alphac}$ there are one-sided syntactic ideals $I^{\ell}_{\alphai}$ and  $I^{r}_{\alphai}$. Here
\begin{equation*}
    I^{\ell}_{\alphai}=\{x\in\klS|\alphai(yx)=0 \ \forall y\in \klS\} \ \ \mathrm{and} \ \  I^{r}_{\alphai}=\{x\in\klS|\alphai(xy)=0 \ \forall y\in \klS\}
\end{equation*} 
are left and right ideals in $\klS$, respectively. 

\begin{figure} 
\begin{center}
\includegraphics[scale=1.0]{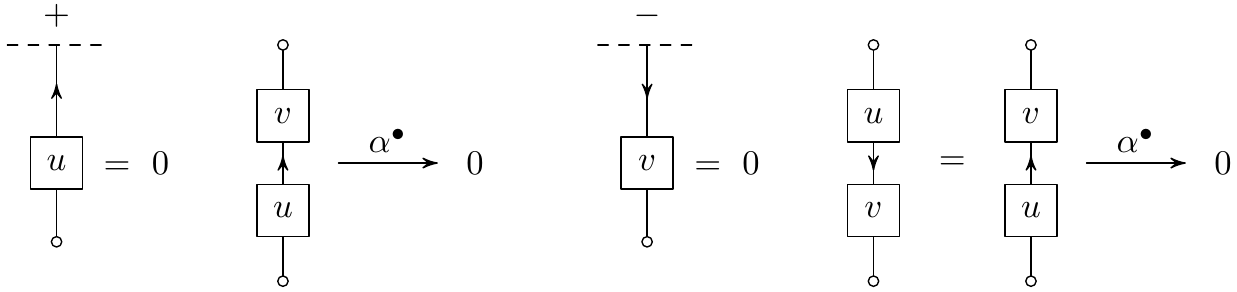}
\caption{Left: Element $\cob^+(u)$ of $\Hom(\emptyset,(+))$ is  set to zero  in $\sCCI_{\alpha}$ for $u\in I^{\ell}_{\alphai}$. This relation is compatible with evaluations of floating  diagrams, since $\alphai(vu)=0$ $\forall v\in\klS$, second left. Right:
Element $\cob^-(v)$ of $\Hom(\emptyset,(-))$ is  set to zero  in $\sCCI_{\alpha}$ for $v\in I^{r}_{\alphai}$. This relation  is also compatible with evaluations  of floating  diagrams, since
$\alphai(vu)=0$ $\forall u\in\klS$ for such  $v$, see the last equation  on the right.}
\label{fig4_11}
\end{center}
\end{figure}

For  $u\in  \klS$ denote by $\cob^+(u)$ the element of $\Hom(\emptyset,(+))$ given by putting $u$ on an interval at its "out" floating endpoint, see Figure~\ref{fig4_11} left. Define $\cob^-(u)$ likewise, see Figure~\ref{fig4_11} right. 
We add relations that $\cob^+(u)=0$ for $u\in I^{\ell}_{\alphai}$ and  $\cob^-(v)=0$ for $v\in I^{r}_{\alphai}$. This finishes our  definition of category $\sCCI_{\alpha}$.

Note  that ideals $I_{\alpha},I^{\ell}_{\alphai},I^{r}_{\alphai}$ have finite codimensions in $\klS$,  and  each of these ideals is finitely-generated. In particular, one can restrict to adding finitely many relations to $\vCCI_{\alpha}$ to get the "skein" category $\sCCI_{\alpha}$. 

Due to consistency of these relations with the evaluation  $\alpha$ on floating components we can describe a basis in the hom spaces in the category $\sCCI_{\alpha}$, as follows. Choose subsets   $B_{\alpha},B^{\ell}_{\alphai},B^r_{\alphai}\subset \klS$ that descend to bases of $A_{\alpha},\klS/I^{\ell}_{\alphai}$ and $I^r_{\alphai}\slash \klS$, respectively. 

Recall the basis $\mcI^S(\undep,\undep')$ of the hom  space from $\undep$ to  $\undep'$  in $\vCCI_{\alpha}$ constructed earlier. It consists of a floating cobordism $x$ from $\undep$ to $\undep'$ with various monomials written on components of the cobordism (all components  are viewable). Define the set $B_{\alpha}(\undep,\undep')$ to also consists of floating cobordisms from  $\undep$ to  $\undep'$, but now we write an  element  of one of the three  sets $B_{\alpha},B^{\ell}_{\alphai},B^r_{\alphai}$ on each component of $c$, depending on its  type: 
\begin{itemize}
    \item If a component has no floating endpoints, thus connects  two boundary points (at the top or bottom  boundary, or both), put an element of $B_{\alpha}$ along  it. 
    \item If a component  has a floating endpoint and is oriented away from this endpoint, put an  element  of $B^{\ell}_{\alphai}$ along this  component. 
    \item If a component  has a floating endpoint and is oriented towards it, put an  element  of $B^{r}_{\alphai}$ along this component.
\end{itemize}
An undecorated viewable cobordism $x$ with $n_1,n_2,n_3$ components of these three types, respectively, admits $|B_{\alpha}|^{n_1}|B^{\ell}_{\alphai}|^{n_2}|B^r_{\alphai}|^{n_3}$ possible decorations. The set $B_{\alpha}(\undep,\undep')$ is the union of these decorated cobordisms, where we start with any viewable undecorated cobordism $x$  from $\undep$ to $\undep'$ and decorate it in all possible such ways. The set $B_{\alpha}(\undep,\undep')$ is finite. 

\begin{figure} 
\begin{center}
\includegraphics[scale=1.0]{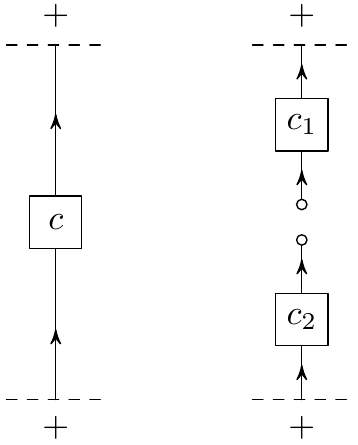}
\caption{There are two  types of cobordisms from $(+)$ to $(+)$ in $\CCI$. Mirroring that decomposition, basis $B_{\alpha}((+),(+))$ consists of elements of $c\in B_{\alpha}$ placed on through strand and pairs of elements $c_1\in B^{\ell}_{\alphai}$ and $c_2\in B^{r}_{\alphai}$ placed on the two strands with floating endpoints.}
\label{fig4_12}
\end{center}
\end{figure}

For example, $B_{\alpha}((+),(+))$ has cardinality $|B_{\alpha}|+|B^{\ell}_{\alphai}| \cdot |B^r_{\alphai}|$ and consists  of diagrams of  two  types, see Figure~\ref{fig4_12}. 

\vspace{0.1in} 

\begin{prop} The set  $B_{\alpha}(\undep,\undep')$ is a basis of the hom  space $\Hom_{\sCCI_{\alpha}}(\undep,\undep')$ in the  category $\sCCI_{\alpha}$. 
\end{prop} 

This construction gives a basis in hom spaces of $\sCCI_{\alpha}$  for  non-recognizable $\alpha$  as well, but then $A_{\alpha}$ and  hom spaces are infinite-dimensional. Recall that we restrict to considering recognizable $\alpha$  for most of this section. 

The endomorphism ring  of $(+)$ in $\sCCI_{\alpha}$ contains a two-sided ideal isomorphic to the  tensor  product  
${\klS/I^{\ell}_{\alphai}\otimes I^r_{\alphai}\backslash \klS}$, with the quotient algebra isomorphic to $A_{\alpha}$, so there is an exact sequence of $\klS$-bimodules 
\begin{equation}
    0 \lra \klS/I^{\ell}_{\alphai}\otimes I^r_{\alphai}\backslash \klS  \lra \End_{\sCCI_{\alpha}}((+)) \lra A_{\alpha} \lra 0 .
\end{equation}
The quotient map onto $A_{\alpha}$ admits a section, and $A_{\alpha}$ is naturally a subalgebra of $\End_{\sCCI_{\alpha}}((+))$. This decomposition corresponds  to two types of endomorphisms of $(+)$ in $\CCI$ (without floating  endpoints versus having two floating endpoints) and corresponding bases in endomorphisms of $(+)$ in $\sCCI_{\alpha}$, see Figure~\ref{fig4_12}. 

\vspace{0.1in} 

(4) {\it Gligible quotient category  $\CCI_{\alpha}$.} 
The \emph{trace} $\tral(x)$ of a cobordism $x$ from $\undep$ to $\undep$ is defined in the same way as for cobordisms in the smaller  category $\CCC_{\alpha}$, by closing $x$ into  a floating cobordism $ \widehat{x}$, see Figure~\ref{fig3_3} and evaluating  via $\alpha$: 
\begin{equation*}
  \tral(x) := \alpha(\widehat{x}). 
\end{equation*} 
This operation  extends to a $\kk$-linear trace on $\kk\CCI_{\alpha}$ that descends to a trace on $\vCCI_{\alpha}$ and $\sCCI_{\alpha}$: 
\begin{equation*}
  \End_{\kk\CCI}(\undep)\lra \End_{\vCCI_{\alpha}}(\undep)\lra  \End_{\sCCI_{\alpha}}(\undep) \stackrel{\tral}{\lra}\kk. 
\end{equation*}

The trace is symmetric. The two-sided ideal $\widetilde{J}_{\alpha}\subset \sCCI_{\alpha}$ of \emph{negligible morphisms} is defined as usual, see Section~\ref{subsec_from} for the definition of negligible ideal in the subcategory $\sCCC_{\alpha}$  of $\sCCI_{\alpha}$.

Define the quotient category  
\begin{equation*}
    \CCI_{\alpha} := \sCCI_{\alpha}/J_{\alpha}.
\end{equation*}
 
 The quotient category $\CCI_{\alpha}$ has finite-dimensional hom spaces, as does  $\sCCI_{\alpha}$, since $\alpha$ is recognizable. The trace is nondegenerate on $\CCI_{\alpha}$ and defines perfect bilinear  pairings
 \begin{equation*}
     \Hom(\undep,\undep') \otimes \Hom(\undep',\undep) \lra \kk
 \end{equation*}
 on its hom spaces. We call $\CCI_{\alpha}$ the  \emph{gligible quotient} of $\sCCI_{\alpha}$, having modded out by  the ideal of negligible morphisms. 
 
\vspace{0.1in}

Up to an isomorphism, the state space 
\begin{equation}
    A_{\alpha}(\undep) \ := \ \Hom_{\CCI_{\alpha}}(\emptyset,\undep)
\end{equation}
depends only on the number of pluses and minuses in $\undep$ and $A_{\alpha}(\undep)\cong A_{\alpha}(+^n-^m)$, where $n$ and $m$ is the number of pluses and minuses in $\undep$. Summing $ A_{\alpha}(+^n-^m)$ over $n,m\ge 0$ one get a bigraded associative algebra with $S_n\times S_m$ action on the homogeneous $(n,m)$ component with the properties similar to that of a \emph{tca} algebra~\cite{SSn}. 

\vspace{0.1in} 

{\it (5) The Deligne category and its gligible quotient.} From the skein category $\sCCI_{\alpha}$ we can pass to its additive Karoubi closure 
\begin{equation*}
    \dCCI_{\alpha} \ := \ \Kar(\sCCI^{\oplus}_{\alpha}),
\end{equation*}
which is the analogue of the Deligne  category. The quotient of $\dCCI_{\alpha}$ by the ideal $\mathcal{D}\widetilde{\mathcal{J}}_{\alpha}$ of negligible morphisms, 
\begin{equation*}
    \udCCI_{\alpha} \ := \ \dCCI_{\alpha}/ \mathcal{D}\widetilde{\mathcal{J}}_{\alpha},
\end{equation*}
is equivalent  to the additive Karoubi closure of $\CCI_{\alpha}$.

\vspace{0.1in} 

{\it Summary:}
To  summarize, the following  categories are assigned to a recognizable pair $\alpha$ as  in (\ref{eq_pair_ser}):  
\begin{itemize}
\item The category $\vCCI_{\alpha}$ of viewable cobordisms with the $\alpha$-evaluation  of floating (or closed) components. 
\item The skein category $\sCCI_{\alpha}$ where closed (floating) $S$-decorated intervals  and circles are evaluated via $\alpha$ and elements of the syntactic ideal $I_{\alpha}$ evaluate to  zero when placed along any interval in a cobordism. Furthemore, elements of left and right syntactic ideals $I_{\alphai}^{\ell}$ and   $I_{\alphai}^{r}$ evaluate to zero when placed at the beginning or end of an interval with the corresponding endpoint floating. 
\item The  quotient $\CCI_{\alpha}$ of $\sCCI_{\alpha}$ by the two-sided ideal of negligible morphisms. We also call $\CCI_{\alpha}$ the \emph{gligible} category or the gligible quotient. Hom spaces in $\CCI_{\alpha}$ come with nondegenerate bilinear forms 
\begin{equation*}
    \Hom(\undep,\undep')\otimes \Hom(\undep',\undep)\lra \kk,
\end{equation*}
where  $\undep,\undep'$ are objects $\CCI_{\alpha}$, sequences of pluses and minuses  describing oriented  zero-manifolds that are the  source and  the target  of decorated  one-cobordisms. The universal construction, in this case, assigns the vector space
$\Hom_{\CCI_{\alpha}}(\emptyset,\undep)$ of morphisms from the empty  zero-manifold $\emptyset$ to $\undep$    to the oriented zero-manifold $\undep$. 
\item Additive Karoubi closure $\dCCI_{\alpha}$ of $\sCCI_{\alpha}$, analogous to  the Deligne category. The quotient of $\dCCI_{\alpha}$ by the ideal of negligible morphisms is denoted $\udCCI_{\alpha}$. 
\end{itemize}

We arrange these  categories and functors, for   recognizable $\alpha=\alphap$,
into  the following diagram, with a commutative square on the right: 

\begin{equation} \label{eq_seq_cd_3}
\begin{CD}
\CCI  @>>> \kk\CCI @>>> \vCCI_{\alpha} @>>> \sCCI_{\alpha} @>>>  \dCCI_{\alpha} \\
@.   @. @.   @VVV   @VVV  \\    
 @.  @.   @.  \CCI_{\alpha}  @>>>  \udCCI_{\,\alpha}
\end{CD}
\end{equation}
Properties of the categories in the analogous diagram (\ref{eq_seq_cd_1})  in Section~\ref{subsec_summ}, as explained in the paragraph following (\ref{eq_seq_cd_1}), hold for the categories in (\ref{eq_seq_cd_3}) as well. 

\vspace{0.1in} 

The category  $\CCI$ and categories built out of  it  require a pair of series $\alpha=\alphap$ for  evaluation. When  working with the subcategory $\CCC$ of cobordisms without floating endpoints,  only circles appear as connected components  of floating cobordisms, and series $\alphac$ is needed  for evaluation. Instead of working with $\CCC$, one can use $\CCI$ but set the  connected component $\alphai=0$, so  that $\alpha=(0,\alphac)$. Then in the viewable category $\vCCI_{\alpha}$ any cobordism that contains a floating interval evaluates to zero. Syntactic ideals $I^{\ell}_{\alphai},I^r_{\alphai}=\klS$, and  in the skein category $\sCCI_{\alpha}$ any cobordism containing an interval  (including  viewable intervals, with  one boundary and  one  floating endpoint) evaluates to $0$. This results in equivalences of categories 
\begin{equation} 
\sCCI_{(0,\alphac)}\cong \sCCC_{\alphac}, \ \ \ 
\CCI_{(0,\alphac)}\cong \CCC_{\alphac}, \ \ \ 
\dCCI_{(0,\alphac)}\cong \dCCC_{\alphac}, \ \ \ 
\udCCI_{(0,\alphac)}\cong \udCCC_{\alphac}, 
\end{equation}
that, furthermore, respect commutative squares of categories in diagrams (\ref{eq_seq_cd_1}) and (\ref{eq_seq_cd_3}).

Thus, the construction in Section~\ref{sec_decorated} of the skein category $\sCCC_{\alpha}$, the gligible quotient category $\CCC_{\alpha}$ and  other categories defined there and   associated to symmetric series $\alpha$ can be considered a  special case of the construction of the present section, specializing  to the pair  $(0,\alpha)$ with the first recognizable series  in the pair being  zero. 

\vspace{0.1in}

Alternatively, one  can set $\alphac$  to zero and consider a pair
$\alpha=(\alphai,0)$. In the viewable cobordism category $\vCCI_{\alpha}$  for this $\alpha$ a circle (necessarily  floating, and with any decoration) evaluates to  $0$. Only the  ideals $I_{\alphai}, I^{\ell}_{\alphai},I^r_{\alphai}$ (two-sided, left  and right, respectively)  are  used in the definition of the skein category $\sCCI_{\alpha}$. A decorated $U$-turn  as   in Figure~\ref{fig4_3}, case (2) or (3), may be non-zero in $\CCI_{\alpha}$, since it may be coupled  to two intervals on the other side with  a non-zero evaluation. 

\vspace{0.1in} 

When $\alpha=(\alphai,0)$, another  approach is to restrict possible cobordisms and disallow $U$-turns as cobordisms. Then  a cobordism $c$ must have no critical points under the natural projection onto the interval $[0,1]$  under which $\partial_i c$ projects onto  $i$, for $i=0,1$. When components of cobordisms  are unoriented, such restricted cobordisms appear in~\cite{KS1} in a categorification of the polynomial ring (without dot decorations) and in~\cite{KT} in a categorification of $\Z[1/2]$ and potential categorifications of $\Z[1/n]$ as monoidal envelopes of certain Leavitt  path algebras and  the "one-sided inverse" algebra $\kk\langle a,b\rangle/(ab-1)$. The latter cobordisms carry dot gecorations, corresponding to the generators of these algebras. 

\vspace{0.1in} 

When $S=\emptyset$, the evaluation again reduces to two numbers (evaluations of  the oriented interval and oriented circle), and the skein category $\sCCI$ is the oriented partial Brauer category, see the remark below. When $S=\{s\}$ has cardinality  one, recognizable series $\alpha$ is encoded by two rational functions $Z_{\alphai}(T)$, $Z_{\alphac}(T)$ in a single variable $T$. 

 \vspace{0.1in}
 
\emph{Remark:} Instead  of power series $\alpha=(\alphai,\alphac)$ in noncommuting variables one can instead start with an associative $\kk$-algebra $B$ and two $\kk$-linear maps 
\begin{equation}
\alphai  \ : \  B\lra \kk, \  \ \alphac \ : \ B \lra \kk
\end{equation} 
such that $\alphac$ is symmetric, $\alphac(ab)=\alphac(ba) $, $a,b\in B$. Two-sided syntactic ideals $I_{\alphai},I_{\alphac}\subset  B$, their intersection $I_{\alpha} := I_{\alphai}\cap I_{\alphac}$, and one-sided syntactic  ideals $I^r_{\alphai}, I^{\ell}_{\alphai}$ are defined in the same way as for noncommutative series. 

The nondegenerate case is that of  $I_{\alpha}=0$ being the zero ideal in $B$, but the arbitrary case can be reduced to it by passing to the quotient $B/I_{\alpha}$. Recognizable case corresponds to finite-dimensional $B$. Analogues of all categories in (\ref{eq_seq_cd_3}) can be defined for such pair of traces on a finite-dimensional $B$. Defects on cobordisms are now labelled by elements of $B$ rather than by elements of $S$. The difference from noncommutative recognizable power series is that one does not pick any particular set  of generators $S$ of $B$, working with the entire $B$ instead, but the resulting  categories, starting with the category $\sCCI_{\alpha}$ in (\ref{eq_seq_cd_3}), are equivalent to the ones built from a noncommutative power series once a set $S$ of generators of $B$ is chosen. 

 \vspace{0.1in}
 
\emph{Remark:} 
It is straightforward to modify the constructions of this section to the case of unoriented one-manifolds with floating endpoints and  $S$-decorated dots. If, in addition, $S=\emptyset$, there are no dots and floating cobordisms reduce to unions  of intervals and circles. The evaluation is then a pair $(\alpha^{\bullet}(1),\alpha^{\circ}(1))$ of elements of $\kk$ and the unoriented skein category $\sCCI_{\alpha}$  is the partial Brauer category~\cite{MM}, also known as the  rook-Brauer category~\cite{HdM}. 



\end{document}